\documentstyle[12pt,leqno]{article}

\newcommand{\sect}[1]{\section{#1}\setcounter{equation}{0}}

\begin{document}

\newtheorem{theorem}{Theorem}[section]
\newtheorem{corollary}[theorem]{Corollary}
\newtheorem{lemma}[theorem]{Lemma}
\newtheorem{proposition}[theorem]{Proposition}
\newtheorem{example}[theorem]{\sl Example}
\newtheorem{remark}[theorem]{\sl Remark}
\newtheorem{definition}[theorem]{\sl Definition}
\newtheorem{algorithm}[theorem]{\sl Algorithm}
\newcommand{\qed}{\ \ \rule{1ex}{1ex}}
\newcommand{\prob}{{\bf P}}  
\newcommand{\E}{{\bf E}}  
\newcommand{\lc}{\left\lceil}
\newcommand{\lf}{\left\lfloor}
\newcommand{\rc}{\right\rceil}
\newcommand{\rf}{\right\rfloor}
\newcommand{\proof}{\par\noindent  {\bf  Proof.  }}
\newcommand{\bd}[1]{{\bf  #1}}
\newcommand{\LLL}{{\bf  \Lambda}}
\newcommand{\OO}{{\bf  \Omega}}
\newcommand{\DD}{{\bf  D}}
\newcommand{\EE}{{\bf  E}}
\newcommand{\KK}{{\bf  K}}
\newcommand{\LL}{{\bf  L}}
\newcommand{\RR}{{\bf  R}}
\newcommand{\NN}{{\bf  N}}
\newcommand{\ZZ}{{\bf  Z}}
\newcommand{\BB}{{\bf  B}}
\newcommand{\CC}{{\bf  C}}
\newcommand{\MM}{{\bf  M}}
\newcommand{\PP}{{\bf  P}}
\newcommand{\QQ}{{\bf  Q}}
\newcommand{\SS}{{\bf  S}}
\newcommand{\TT}{{\bf  T}}
\newcommand{\UU}{{\bf  U}}
\newcommand{\VV}{{\bf  V}}
\newcommand{\XX}{{\bf  X}}
\newcommand{\XXp}{{\bf X'}}
\newcommand{\YY}{{\bf  Y}}
\newcommand{\WW}{{\bf  W}}
\newcommand{\AAA}{{\bf  A}}
\newcommand{\pp}{{\bf  p}}
\newcommand{\uu}{{\bf  u}}
\newcommand{\xx}{{\bf  x}}
\newcommand{\yy}{{\bf  y}}
\newcommand{\zz}{{\bf  z}}
\newcommand{\ww}{{\bf  w}}
\newcommand{\Sj}{{{\cal S}_j}}
\newcommand{\Kt}{{\widetilde{K}}}
\newcommand{\Lt}{{\widetilde{L}}}
\newcommand{\Pt}{{\widetilde{P}}}
\newcommand{\PPt}{{\widetilde{\PP}}}
\newcommand{\MMt}{{\widetilde{\MM}}}
\newcommand{\Var}{{\rm Var}}
\newcommand{\diag}{{\rm diag}}
\newcommand{\Exp}{{\rm Exp}}
\newcommand{\ave}{{\rm E}}
\newcommand{\id}{{\rm id}}
\newcommand{\rev}{{\rm rev}}
\newcommand{\Xt}{{\widetilde{X}}}
\newcommand{\Yt}{{\widetilde{Y}}}
\newcommand{\XXt}{{\widetilde{\XX}}}
\newcommand{\YYt}{{\widetilde{\YY}}}
\newcommand{\SSS}{{S}} 
\newcommand{\gxy}{g_{xy}}
\newcommand{\Gamxy}{{\Gamma(x, y)}}
\newcommand{\vecXX}{{\vec{\XX}}}
\newcommand{\vecUU}{{\vec{\UU}}}
\newcommand{\vecVV}{{\vec{\VV}}}
\newcommand{\vecWW}{{\vec{\WW}}}
\newcommand{\vecYY}{{\vec{\YY}}}
\newcommand{\vecDD}{{\vec{\DD}}}
\newcommand{\vecD}{{\vec{D}}}
\newcommand{\vecL}{{\vec{L}}}
\newcommand{\vecV}{{\vec{V}}}
\newcommand{\vecU}{{\vec{U}}}
\newcommand{\vecY}{{\vec{Y}}}
\newcommand{\vecX}{{\vec{X}}}
\newcommand{\vecw}{{\vec{w}}}
\newcommand{\vecx}{{\vec{x}}}
\newcommand{\vecy}{{\vec{y}}}
\newcommand{\vecu}{{\vec{u}}}
\newcommand{\vece}{{\vec{e}}}
\newcommand{\vecep}{{\vece\,'}}
\newcommand{\Qe}{{Q(\vece\,)}}
\newcommand{\Qep}{{Q(\vecep)}}
\newcommand{\we}{{w(\vece\,)}}
\newcommand{\wep}{{w(\vecep)}}
\newcommand{\Kwe}{{K_w(\vece\,)}}
\newcommand{\phie}{{\phi(\vece\,)}}
\newcommand{\psie}{{\psi(\vece\,)}}
\newcommand{\be}{{b(\vece\,)}}
\newcommand{\vecE}{{\vec{E}}}
\newcommand{\GamxyQ}{{|\Gamma(x, y)|_Q}}
\newcommand{\Gamxyw}{{|\Gamma(x, y)|_w}}
\newcommand{\gamstar}{{\gamma_*}}
\newcommand{\gamstart}{{\gamma^t_*}}
\newcommand{\gamstarone}{{\gamma^1_*}}
\newcommand{\gamstartwo}{{\gamma^2_*}}
\newcommand{\gamstarthree}{{\gamma^3_*}}
\newcommand{\gamstarfour}{{\gamma^4_*}}
\newcommand{\bte}{{b^t(\vece\,)}}
\newcommand{\bonee}{{b^1(\vece\,)}}
\newcommand{\btwoe}{{b^2(\vece\,)}}
\newcommand{\bthreee}{{b^3(\vece\,)}}
\newcommand{\bfoure}{{b^4(\vece\,)}}
\newcommand{\bthreeae}{{b^{3\mbox{\scriptsize a}}(\vece\,)}}
\newcommand{\bthreebe}{{b^{3\mbox{\scriptsize b}}(\vece\,)}}
\newcommand{\bfourae}{{b^{4\mbox{\scriptsize a}}(\vece\,)}}
\newcommand{\bfourbe}{{b^{4\mbox{\scriptsize b}}(\vece\,)}}
\newcommand{\Ac}{{\cal A}}
\newcommand{\Bc}{{\cal B}}
\newcommand{\Cc}{{\cal C}}
\newcommand{\Dc}{{\cal D}}
\newcommand{\Ec}{{\cal E}}
\newcommand{\Fc}{{\cal F}}
\newcommand{\Gc}{{\cal G}}
\newcommand{\Lc}{{\cal L}}
\newcommand{\Nc}{{\cal N}}
\newcommand{\Sc}{{\cal S}}
\newcommand{\Uc}{{\cal U}}
\newcommand{\Xc}{{\cal X}}
\newcommand{\Zc}{{\cal Z}}
\newcommand{\Ft}{\Fc_t}
\newcommand{\Fi}{\Fc_{\infty}}
\newcommand{\Fn}{\Fc_n}
\newcommand{\Tt}{{\widetilde{T}}}
\newcommand{\ft}{{\tilde{f}}}
\newcommand{\Jt}{{\widetilde{J}}}
\newcommand{\tp}{t+}
\newcommand{\tm}{t-}
\newcommand{\vdn}{\| \pi_n - \pi \|}
\newcommand{\vdt}{\| \pi_t - \pi \|}
\newcommand{\sep}{{\rm sep}}
\newcommand{\Fets}{\Fc^*_{=t}}
\newcommand{\ppi}{{\mbox{\boldmath $\pi$}}}
\newcommand{\ssigma}{{\mbox{\boldmath $\sigma$}}}
\newcommand{\PPh}{\PP^{(h)}}
\newcommand{\GG}{{\bf  G}}
\newcommand{\II}{{\bf  I}}
\newcommand{\Qsp}{{Q^*}'}
\newcommand{\sn}{\sigma_n}
\newcommand{\tn}{\tau_n}
\newcommand{\tnm}{\tau_{n-1}}
\newcommand{\gt}{{\tilde{g}}}
\newcommand{\Gt}{{\widetilde{G}}}
\newcommand{\pimin}{\pi_{\min}}
\newcommand{\zh}{\hat{0}}
\newcommand{\oh}{\hat{1}}
\newcommand{\Th}{\hat{T}}
\newcommand{\th}{\hat{t}}
\newcommand{\pih}{\hat{\pi}}
\newcommand{\taut}{\tilde{\tau}}
\newcommand{\hti}{\tilde{h}}
\newcommand{\bege}{\begin{equation}}
\newcommand{\ene}{\end{equation}}
\newcommand{\begp}{\begin{proposition}}
\newcommand{\enp}{\end{proposition}}
\newcommand{\begt}{\begin{theorem}}
\newcommand{\ent}{\end{theorem}}
\newcommand{\begl}{\begin{lemma}}
\newcommand{\enl}{\end{lemma}}
\newcommand{\begc}{\begin{corollary}}
\newcommand{\enc}{\end{corollary}}
\newcommand{\begr}{\begin{remark}}
\newcommand{\enr}{\end{remark}}
\newcommand{\begd}{\begin{definition}}
\newcommand{\enf}{\end{definition}}
\newcommand{\begx}{\begin{example}}
\newcommand{\enx}{\end{example}}
\newcommand{\bega}{\begin{array}}
\newcommand{\ena}{\end{array}}
\newcommand{\Line}{$\underline{\qquad\qquad\qquad\qquad}$}
\newcommand{\X}{$\qquad$}
\newcommand{\lb}{\left\{}
\newcommand{\rb}{\right\}}
\newcommand{\lsb}{\left[}
\newcommand{\rsb}{\right]}
\newcommand{\lp}{\left(}
\newcommand{\rp}{\right)}
\newcommand{\ls}{\left|}
\newcommand{\rs}{\right|}
\newcommand{\lss}{\left\|}
\newcommand{\rss}{\right\|}
\newcommand{\vs}{\vspace{\smallskipamount}\noindent}
\newcommand{\vm}{\vspace{\medskipamount}\noindent}
\newcommand{\vb}{\vspace{\bigskipamount}\noindent}
\newcommand{\ra}{\rightarrow}
\newcommand{\la}{\leftarrow}
\newcommand{\implies}{\Longrightarrow}
\newcommand{\s}{\hspace*{5ex}}
\newcommand{\While}{{\tt while\ }}
\newcommand{\Repeat}{{\tt repeat\ }}
\newcommand{\Until}{{\tt until\ }}
\newcommand{\For}{{\tt for\ }}
\newcommand{\To}{{\tt \ to\ }}
\newcommand{\If}{{\tt if\ }}
\newcommand{\Then}{{\tt then}}
\newcommand{\Else}{{\tt else}}
\newcommand{\Return}{{\tt return\ }}
\newcommand{\Failure}{{\tt Failure}}
\newcommand{\set}{\leftarrow}
\newcommand{\unknown}{{\tt UNKNOWN}}
\newcommand{\defined}{{\tt DEFINED}}
\newcommand{\rmap}{\mbox{\tt RandomMap()}}
\newcommand{\next}{\mbox{\tt NextState()}}
\newcommand{\curr}{{\tt CurrentState()}}
\newcommand{\oF}{{\overline F}}
\newcommand{\Tmix}{T_{\text{mix}}}
\newenvironment{code}{\begin{list}{\hspace*{0em}}{\renewcommand{\parsep}{-2pt}}}{\end{list}}
 
\newcommand{\sfrac}[2]{{\textstyle\frac{#1}{#2}}}

\setcounter{page}{-1}

\author{
{\bf James Allen Fill}
$^*$
and {\bf Motoya Machida}\thanks{
These authors' research has been supported in part by NSF grants DMS--9626756 and
DMS--9803780, and by the Acheson J.~Duncan Fund for the Advancement of Research in
Statistics.}
\\
Department of Mathematical Sciences \\
The Johns Hopkins University \\
Baltimore, MD 21218--2682 \\
U.~S.~A. \\
{\tt jimfill@jhu.edu} \\
{\tt machida@mts.jhu.edu} \\
{\ }
\and {\bf Duncan J.~Murdoch}\thanks{
These authors' research has been supported in part by NSERC of Canada.}
\\
Department of Statistical and Actuarial Sciences \\
University of Western Ontario \\
London, Ontario N6G 2E9 \\
Canada \\
{\tt murdoch@fisher.stats.uwo.ca} \\
{\ }
\and {\bf Jeffrey S.~Rosenthal}$^{\dag}$
\\
Department of Statistics \\
University of Toronto \\
Toronto, Ontario M5S 3G3 \\
Canada \\
{\tt jeff@math.toronto.edu}
}
\title{{\bf Extension of Fill's perfect rejection sampling algorithm to general
chains}}
\date{April, 1999; last revised July~20, 2000}
\maketitle

\thispagestyle{empty}

\newpage

\thispagestyle{empty}

\begin{center}
{\bf Abstract}
\end{center}
 
By developing and applying a broad framework for rejection sampling using
auxiliary randomness, we provide an extension of the perfect sampling algorithm of
Fill (1998) to general chains on quite general state spaces, and describe how use
of bounding processes can ease computational burden.  Along the way, we unearth a
simple connection between the Coupling From The Past (CFTP) algorithm originated by
Propp and Wilson (1996) and our extension of Fill's algorithm.
\bigskip
 
\bigskip

\addtolength{\topmargin}{+0.5in}

\newpage
 
\sect{Introduction}
\label{intro}

\hspace{\parindent}
Markov chain Monte Carlo (MCMC) methods have become extremely popular
for Bayesian inference problems (see e.g.\ Gelfand and Smith~\cite{GS}, Smith
and Roberts~\cite{SR}, Tierney~\cite{Tierney}, Gilks et al.~\cite{GRS}), and for
problems in other areas, such as spatial statistics, statistical physics, and
computer science (see e.g.\ Fill~\cite{Fill} or Propp and Wilson~\cite{PWexact}
for pointers to the literature) as a way of sampling approximately from a
complicated unknown probability distribution~$\pi$.  
An MCMC algorithm constructs
a Markov chain with one-step transition kernel~$K$ and stationary
distribution~$\pi$; if the chain is run long enough, then under reasonably weak
conditions (cf.\ Tierney~\cite{Tierney}) it will converge in distribution 
to~$\pi$, facilitating approximate sampling.

One difficulty with these methods is that it is difficult to 
assess convergence to
stationarity.  This necessitates the use of difficult theoretical analysis
(e.g.,\ Meyn and Tweedie~\cite{MT}, Rosenthal~\cite{Rosenthal}) or problematic
convergence diagnostics (Cowles and Carlin~\cite{CC}, Brooks, et al.~\cite{BDR})
to draw reliable samples and do proper inference.

An interesting alternative algorithm, called {\em coupling from the
past\/} (CFTP), was introduced by Propp and Wilson~\cite{PWexact}
(see also~\cite{PWuser} and~\cite{PWhowto}) and has been studied
and used by a number of authors (including Kendall~\cite{Kendall},
M{\o}ller~\cite{Moller}, Murdoch and Green~\cite{MG}, Foss and
Tweedie~\cite{FT}, Kendall and Th\"onnes~\cite{KT}, Corcoran
and Tweedie~\cite{CT}, Kendall and M{\o}ller~\cite{KM}, Green and
Murdoch~\cite{GM}, and Murdoch and Rosenthal~\cite{MR}).  By searching
backwards in time until paths from all starting states have coalesced,
this algorithm uses the Markov kernel~$K$ to sample {\em exactly\/}
from~$\pi$.

Another method of perfect simulation, for finite-state stochastically
monotone chains, was proposed by Fill~\cite{Fill}.  Fill's algorithm
is a form of rejection sampling.  This algorithm was later extended
by M{\o}ller and Schladitz~\cite{MS} and Th\"onnes~\cite{Thonnes}
to non-finite chains, motivated by applications to spatial point
processes.  Fill's algorithm has the advantage over CFTP of removing
the correlation between the length of the run and the returned value,
which eliminates bias introduced by an impatient user or a system
crash and so is ``interruptible''. However, it has been used only for
stochastically monotone chains, making heavy use of the ordering of
state space elements.  In his paper, Fill~\cite{Fill} indicated that
his algorithm could be suitably modified to allow for the treatment of
``anti-monotone'' chains and (see his Section 11.2) indeed to generic
chains.  A valuable background resource on perfect sampling methods is
the annotated bibliography maintained by Wilson~\cite{Wilson}.

The goal of the present paper is to discuss the modifications to
Fill's algorithm needed to apply it to generic chains, on general (not
necessarily finite) state spaces.  Our basic algorithm is presented
in Section~\ref{alginbrief} as Algorithm~\ref{nontech}.
An infinite-time-window version of
Algorithm~\ref{nontech} (namely, Algorithm~\ref{altalg}) is presented
in Section~\ref{altalgsect}.
A simple illustrative example is presented in Section~\ref{toy}, while
rigorous mathematical details are
presented in Section~\ref{frame}.

In Section~\ref{detection} we discuss how the computational burden of
tracking all of the trajectories in Algorithm~\ref{nontech}
can be eased by the use of
coalescence detection events in general and bounding processes in
particular; these processes take on a very simple form (see
Section~\ref{mono}) when the state space is partially ordered and the
transition rule employed is monotone.
A weaker form of monotonicity is also handled in Section~\ref{app3}. 

In Section~\ref{taleof2} we compare Algorithm~\ref{nontech} and CFTP.  We
also present a simple connection between CFTP and Algorithm~\ref{altalg}.
Finally, in Section~\ref{n} we discuss the perfect generation of samples
from~$\pi$ of size larger than~$1$.

We hope that our extension of Fill's algorithm will stimulate further
research into this less-used alternative for perfect MCMC simulation.
\medskip

{\em Notational note:\/} Throughout the paper, we adopt the probabilist's
usual shorthand of writing $\{X \in B\}$ for the event
$\{\omega \in \Omega:\,X(\omega) \in B\}$
when~$X$ is a random element defined on a sample space~$\Omega$.

\sect{The algorithm in brief}
\label{alginbrief}

\hspace{\parindent}
We assume here that our Markov chain may be written in the {\em stochastic
recursive sequence} form
\begin{equation}
\label{SRSeq}
\XX_{s} = \phi(\XX_{s-1},\UU_{s})
\end{equation}
where $(\UU_s)$ is an i.i.d.\ sequence having distribution (say)~$\mu$.  

Omitting technical details, our interruptible algorithm for generic chains is
conceptually quite simple and proceeds as follows.

\begin{algorithm}
\label{nontech}
{\em
Choose and fix a positive integer~$t$, choose an initial state $\XX_t$, and
perform the following routine.  Run the time-reversed chain~$\Kt$
for~$t$ steps [see~(\ref{reversal}) for the formal definition of~$\Kt$],
obtaining $\XX_{t - 1}, \ldots, \XX_0$ in succession.  Then
(conditionally given $\XX_0, \ldots, \XX_t$) 
generate $\UU_1, \ldots, \UU_{t}$ independently,
with~$\UU_s$ chosen from its conditional distribution given~(\ref{SRSeq})
($s = 1, \ldots, t$).
Then, for each element~$x$ of the state space~$\Xc$, compute chains
$(\YY_0(x), \ldots, \YY_t(x))$, with 
$\YY_0(x) := x$ and $\YY_s(x) :=  \phi(\YY_{s-1}(x),\UU_{s})$
for $s=1,2,\ldots,t$.
Note that $\YY_s(\XX_0) = \XX_s$ for $s=1,\ldots,t$.  
Finally, check
whether all the values $\YY_t(x)$, $x \in \Xc$, agree (in which case,
of course, they all equal~$\XX_t$).  If they do, we call this {\em
coalescence\/}, and the routine succeeds and reports
$\WW := \XX_0$ as an observation
from~$\pi$.  If not, then the
routine fails; we then start the routine again with an independent
simulation (perhaps with a fresh choice of~$t$ and~$\XX_t$), and repeat
until the algorithm succeeds.
}
\end{algorithm}

\begin{remark} 
{\em
The algorithm works for $\pi$-almost every deterministic choice of 
initial state $\XX_t$.  Alternatively, the algorithm works provided
one chooses~$\XX_t$
from any distribution absolutely continuous with respect to~$\pi$.
See also Remark~\ref{genremark}(c).
}
\end{remark}

Here is the basic idea of why the algorithm works correctly.  Imagine
(1)~{\em starting\/} the construction with $\XX_0 \sim \pi$ and
independently (2)~simulating $\UU_1, \ldots, \UU_{t}$.  Determination
of coalescence and the value of the coalesced paths at time~$t$ each
rely only on the second piece of randomness.   It follows that,
conditionally given coalescence, $\XX_0$ and $\XX_t$ are
independent.  Hence, conditionally given coalescence and $\XX_t$, we
will still have $\XX_0 \sim \pi$, as desired.  The algorithm constructs
the random variables in a different order, but conditional on
coalescence and the value of $\XX_t$, the joint distributions are the same.

\begin{remark}
\label{nontechremark}
{\em (a)~Note that no assumption is made in Algorithm~\ref{nontech} concerning
monotonicity or discreteness of the state space.

(b) This algorithm is, like Fill's original algorithm~\cite{Fill}, a form
of rejection sampling (see, e.g.,\ Devroye~\cite{Devroye}).
This is explained in Section~2 of~\cite{FMMRshort}.

(c) We have reversed the direction of time, and the roles of the
kernels~$K$ and~$\Kt$, compared to Fill~\cite{Fill}.

(d) Algorithm~\ref{nontech} is interruptible, in the
sense of Fill~\cite{Fill}.

(e) Fill's original algorithm~\cite{Fill}
also incorporated a search for a good value
of~$t$ by doubling the previous value of~$t$ until the first success.
For the most part, we shall not address such issues, instead leaving
the choice of~$t$ entirely up to the user; but see
Section~\ref{altalgsect}.  
} 
\end{remark}

In Section~\ref{frame}, we will carefully discuss the underlying
assumptions for Algorithm~\ref{nontech} and the details of its
implementation, and also establish rigorously that the algorithm works
as desired.  This will be done by first developing, and then applying,
results in a rather more general framework.

\sect{A modified algorithm which searches for $t$}
\label{altalgsect}

\hspace\parindent
Thus far we have been somewhat sketchy about the choice(s) of~$t$ in
Algorithm~\ref{nontech}.  As discussed in Remark~\ref{nontechremark}(e),
one possibility is
to run the repetitions of the basic routine independently, doubling~$t$ at each
stage.  However, another possibility is to continue back in time, reusing the
already imputed values~$\UU_s$ and checking again for coalescence.  (There is an
oblique reference to this alternative in Remark~9.3 of Fill~\cite{Fill}.)  This
idea leads to the following algorithm.

\begin{algorithm}
\label{altalg}
{\em Choose an initial state $\XX_0 \sim \pih$, where $\pih$ is absolutely
continuous with respect to~$\pi$.  Run the time-reversed chain~$\Kt$, obtaining
$\XX_0, \XX_{-1}, \ldots$ in succession.  Then
(conditionally given $\XX_0, \ldots, \XX_t$)
generate $\UU_0, \UU_{-1}, \ldots$ independently,
with~$\UU_s$ chosen from its conditional distribution given~(\ref{SRSeq})
($s = 0, -1, \ldots$).
For $t = 0, 1,
\ldots$ and $x \in \Xc$, set $\YY^{(-t)}_{-t}(x) := x$ and, inductively,
$$
\YY^{(-t)}_s(x) := \phi(\YY^{(-t)}_{s - 1}(x), \UU_s),\ \ \ -t + 1 \leq s \leq 0.
$$
If $\TT < \infty$ is the smallest~$t$ such that
\begin{equation}
\label{altcoalescence}
\mbox{$\YY^{(-t)}_0(x)$, $x \in \Xc$, all agree\ \ \ \ \ \ \ (and
hence all equal $\XX_0$),}
\end{equation}
then the algorithm succeeds and reports~$\WW := \XX_{-\TT}$ as an observation
from~$\pi$.  If there is no such $\TT$, then the algorithm fails.
}
\end{algorithm}

Here is the basic idea of why the algorithm works correctly.  Imagine
(1)~starting the construction with $\WW \sim \pi$, and,
independently, (2)~simulating $\UU_0, \UU_{-1}, \ldots$
[and then, after determining~$\TT$,
setting $\XX_{-\TT}:= \WW$
and $\XX_s := \phi(\XX_{s-1}, \UU_s)$ for $s=-\TT+1,\ldots,0$].  
Determination of~$\TT$ and the value of~$\XX_0$ each
rely only on the second piece of randomness.   It follows that,
conditionally given coalescence, $\XX_{-\TT}$ and $\XX_0$ are
independent.  Hence, conditionally given coalescence and $\XX_0$, we
will still have $\XX_{-\TT} \sim \pi$, as desired.  As before,
the algorithm constructs
the random variables in a different order, but conditional on
coalescence and the value of $\XX_t$, the joint distributions are the same.

\begin{remark}
\label{altremark}
{\em (a)~We need only generate $\XX_0, \XX_{-1}, \ldots, \XX_{-t}$ 
and then impute
$\UU_0$, $\UU_{-1}, \ldots, \UU_{-t + 1}$ in order to check
whether or not~(\ref{altcoalescence}) holds.  Thus, as long
as $\TT < \infty$, the algorithm
runs in finite time.

(b)~One can formulate the algorithm rigorously
in the fashion of
Section~\ref{app1}, and verify that it works properly.  We omit the details.

(c)~Algorithm~\ref{altalg}
is also interruptible:\ \ specifically, $\TT$ and~$\WW$ are 
conditionally independent given success.

(d)~See also the discussion of a ``doubling'' search strategy in
Section~\ref{detection} below.
}
\end{remark}

\sect{A simple illustrative example}
\label{toy}

\hspace\parindent
We illustrate Algorithm~\ref{nontech} for a very simple example
(for which direct sampling from~$\pi$ would be elementary, of course)
and two different
choices of transition rule.  Consider the discrete state space $\Xc = \{0,
1, 2\}$, and let~$\pi$ be uniform on~$\Xc$.  Let~$K$ correspond to
simple symmetric random walk with holding probability~$1/2$ at the endpoints; 
that is, putting $k(x, y) := K(x, \{y\})$,
\begin{eqnarray*}
&& k(0, 0) = k(0, 1) = k(1, 0) = k(1, 2) = k(2, 1) = k(2, 2) = 1/2, \\
&& k(0, 2) = k(1, 1) = k(2, 0) = 0.
\end{eqnarray*}
The stationary distribution is~$\pi$.  As for any ergodic birth-and-death chain,
$K$ is reversible with respect to~$\pi$, i.e.,\ $\Kt = K$.  Before starting the
algorithm, choose a transition rule; this is discussed further below.

For utter simplicity of description, we choose~$t = 2$ and (deterministically)
$\XX_t = 0$ (say); as discussed near the end of Section~\ref{rigor}, a
deterministic start is permissible here.  We then choose $\XX_1 \sim K(0, \cdot)$
and $\XX_0\,|\,\XX_1 \sim K(\XX_1, \cdot)$.  How we proceed from this
juncture depends on what we chose (in advance) for~$\phi$.

One choice is the independent-transitions rule [discussed further in
Remarks \ref{phiremark}(c) and~\ref{subtle}(b) below].  The 
algorithm's routine can then
be run using~$6$ independent random bits: these decide~$\XX_1$ (given~$\XX_2$),
$\XX_0$ (given~$\XX_1$), and the~$4$ transitions in the second (forward)
phase of the routine not already determined from the rule
$$
\XX_{s - 1} \mapsto \XX_s \mbox{\ from time $s - 1$ to time~$s$ ($s = 1, 2$).}
$$
There are thus a total of~$2^6 = 64$ possible
overall simulation results, each having probability~$1/64$.  We check that
exactly~$12$ of these produce coalescence.  Of these~$12$ accepted results,
exactly~$4$ have $\XX_0 = 0$, another~$4$ have $\XX_0 = 1$, and a final~$4$
have $\XX_0 = 2$.  Thus $\PP(\CC) = 12/64 = 3/16$, and we confirm that
$\Lc(\XX_0|C) = \pi$, so the algorithm is working correctly.  (An
identical result holds if we had instead chosen $\XX_t = 1$ or $\XX_t = 2$.)

An alternative choice adapts Remarks~\ref{phiremark}(b) and~\ref{subtle}(c)
to the discrete setting of our present example.
We set
$$
\phi(\cdot, u) = \left\{
 \begin{array}{lll}
  \mbox{the mapping taking\ \ $0,1,2$ to} & \mbox{$0,0,1$,} & \mbox{respectively,
                                                                   if $u = 0$}
\vspace{.1in} \\
  \mbox{the mapping taking\ \ $0,1,2$ to} & \mbox{$1,2,2$,} & \mbox{respectively,
                                                                   if $u = 1$}
 \end{array}
 \right.
$$
where $u$ is uniform on $\{0,1\}$.
Choosing $t = 2$ and $\XX_t = 0$ as before, the algorithm can now be run with
just~$2$ random bits.  In this case we check that exactly~$3$ of the~$4$ possible
simulation results produce coalescence, $1$ each yielding $\XX_0 = 0,1,2$.  Note
that $\PP(\CC) = 3/4$ is much larger for this choice of~$\phi$.  In fact,
since~$\phi$ is a monotone transition rule [see Definition~4.2 in Fill~\cite{Fill}
or (\ref{mondef}) below], for the choice $\XX_t = 0$ it gives the highest
possible value of~$\PP(\CC)$ among all choices of~$\phi$: see Remark~9.3(e) in
Fill~\cite{Fill}.  It also is a best choice when $\XX_t = 2$.  [On a minor negative
note, we observe that $\PP(\CC) = 0$ for the choice $\XX_t = 1$.  Also note that
the $\pi$-average of the acceptance probabilities~$(3/4, 0,
3/4)$, namely, $1/2$, is the probability that forward coupling (or CFTP) done with
the same transition rule gives coalescence within~$2$ time units; this 
corroborates Remark~\ref{genremark}(c) below.]

\begin{remark}
\label{toyrate}
{\em Both choices of~$\phi$ are easily extended to handle simple symmetric random
walk on $\{0, \ldots, n\}$ for any~$n$.  If $\XX_t=0$, then
the second (monotone) choice is again
best possible.  For fixed $c \in (0, \infty)$
and large~$n$, results in Fill~\cite{Fill} and Section~4 of Diaconis and
Fill~\cite{DFAP} imply that, for $t = c n^2$ and $\XX_t=0$, the 
routine's success probability is
approximately~$p(c)$; here~$p(c)$ increases smoothly from~$0$ to~$1$ as $c$
increases from~$0$ to~$\infty$.  We have not attempted the corresponding
asymptotic analysis for the independent-transitions rule.
}
\end{remark}

\sect{Conservative detection of coalescence}
\label{detection}

\hspace{\parindent}
({\em Note:\/} In order to focus on the main ideas, in this
Section~\ref{detection}---as in the previous sections---we
will suppress measure-theoretic and other technical details.
In an earlier draft we provided these details; we trust that the interested reader
will be able to do the same, using the rigorous treatment of Algorithm~\ref{nontech}
in Section~\ref{frame} as a guide.
Also note that the terminology used in this section---detection process, bounding process,
etc.---has varied somewhat in the perfect sampling literature.)
\vspace{0.01in}

Even for large finite state spaces~$\Xc$, determining exactly whether
or not coalescence occurs in Algorithm~\ref{nontech} can be
prohibitively expensive computationally; indeed, in principle this
requires tracking each of the trajectories $\vecYY(x) := (\YY_0(x), \ldots, \YY_t(x))$,
$x \in \Xc$, to completion.  

However, suppose that $\EE$ is an event which is a {\em subset} of 
the coalescence event.  That is, whenever $\EE$ occurs, coalescence occurs
(but perhaps not conversely).  Assume further that $\EE$, like the
coalescence event, is an event whose occurrence (or not) is 
determined solely by $\UU_1,\ldots,\UU_t$.
Then Algorithm~\ref{nontech} remains valid if, instead of accepting 
$\WW$ whenever coalescence occurs, we accept only when the event $\EE$ occurs.
Indeed, the explanation for why Algorithm~\ref{nontech} works goes
through without change in this case.

It follows that, when implementing
Algorithm~\ref{nontech}, it is permissible to use {\em conservative
detection of coalescence\/}, i.e.,\ to accept~$\WW$ as an observation from~$\pi$
if and only if some event~$\EE$ occurs, provided that~$\EE$ is a subset of
the coalescence event and is a function of $\UU_1, \ldots, \UU_t$ only.
We call such an event~$\EE$ a {\em coalescence detection event\/}.

Similar considerations apply to Algorithm~\ref{altalg}.
Indeed, let~$\TT$ be as in that algorithm.  Let $\TT'$ be any other
positive-integer-valued random variable, which is completely determined
by $\UU_0, \UU_{-1}, \ldots$ and is such that $\TT' \ge \TT$.  Then
Algorithm~\ref{altalg} remains valid if we replace $\TT$ by $\TT'$,
i.e.,\ if we report $\WW' := \XX_{-\TT'}$ instead of reporting
$\XX_{-\TT}$.
Indeed, the explanation for why Algorithm~\ref{altalg} works goes
through without change in this case.

For example, we might choose in Algorithm~\ref{altalg} to
let~$\TT'$ be the smallest~$t$ which is a power
of~$2$ such that~(\ref{altcoalescence}) holds and report $\XX_{-\TT'}$ instead.
This has the computational efficiency advantage (similar to that of CFTP)
that we can use a ``doubling'' search strategy, trying
$t=1, 2, 4, 8, \ldots$ in succession,
until we find the first time $t\,( = \TT')$ such that~(\ref{altcoalescence}) holds.

In such a case, and in other cases discussed below,
conservative coalescence detection will often lead to easier and faster application
of our algorithms.  Thus, such detection may be of considerable help in
practical applications.

\subsection{Detection processes}
\label{detectproc}

\hspace{\parindent}
In practice, a coalescence detection event is constructed in terms of
a {\em detection process\/}.  What we mean by this is a stochastic process $\vecDD =
(\DD_0, \ldots, \DD_t)$, defined on the same probability space
as~$\vecUU = (\UU_1, \ldots, \UU_t)$ and~$\vecXX = (\XX_0, \ldots, \XX_t)$, together with a
subset~$\Delta$ of its state space~$\Dc$, such that
\begin{enumerate}
\item[(a)] $\vecDD$ is constructed solely from~$\vecUU$, and
\item[(b)] $\{\DD_s \in \Delta\mbox{\ for some $s \leq t$}\} \subseteq
\{\mbox{$\YY_t(x)$ does not depend on~$x$}\}$.
\end{enumerate}
Then $\EE := \{\DD_s \in \Delta\mbox{\ for some $s \leq t$}\}$ is a coalescence
detection event.

\begin{remark}
\label{detectionremark}
{\em
In practice,
$\vecDD$ usually evolves Markovianly using~$\vecUU$; more precisely, it is
typically the case that there exists deterministic $d_0 \in \Dc$
such that $\DD_0 = d_0$
and~$\vecDD$ has the stochastic recursive sequence form
[paralleling~(\ref{SRSeq})]
$$
\DD_s := \delta(\DD_{s-1}, \UU_s),\ \ \ 1 \leq s \leq t.
$$
}
\end{remark} 

The important consequence is that, having determined the trajectory~$\vecXX$ and
the imputed~$\vecUU$, the user need only follow a single trajectory in the forward
phase of the routine, namely, that of~$\vecDD$.

\begin{example}
\label{MTFandtrees}
{\em
We sketch two illustrative examples of the use of detection processes that do
not immediately fall (see Remark~\ref{detvsbd} below)
into the more specific settings of Sections~\ref{bounding}
or~\ref{mono}.
We hasten to point out, however, that because of the highly special
structure of these two examples, efficient implementation of
Algorithm~\ref{nontech} avoids the use of the forward phase altogether; this is
discussed for example~(a) in Fill~\cite{MTFalgo}.

(a)~Our first example is provided by the move-to-front (MTF) rule studied
in~\cite{MTFalgo}.  Let~$K$ be the Markov kernel corresponding to~MTF with
independent and identically distributed record requests corresponding to
probability weight vector $(w_1, \ldots, w_n)$; see~(2.1) of~\cite{MTFalgo} for
specifics.  The arguments of Section~4 of~\cite{MTFalgo} show that if~$\DD_s$ is
taken to be the set of all records requested at least once among the first~$s$
requests and~$\Delta$ is taken to consist of all $(n - 1)$-element subsets of the
records $1, \ldots, n$, then~$\vecDD$ is a detection process.  Similar detection
processes can be built for the following generalizations of MTF:\ move-to-root for
binary search trees (see Dobrow and Fill~\cite{DF1}~\cite{DF2}) and MTF-like
shuffles of hyperplane arrangement chambers and more general structures (see
Bidigare, et al.~\cite{BHR} and Brown and Diaconis~\cite{BD}).

(b)~A second example of quite similar spirit is provided by the (now well-known)
Markov chain~$(\XX_t)$ for generating a random spanning arborescence of the
underlying weighted directed graph, with vertex set~$\Uc$, of a Markov chain
$(\UU_t)$ with state space~$\Uc$ and kernel~$q$.  Consult
Propp and Wilson~\cite{PWhowto} (who also discuss a more efficient
``cycle-popping'' algorithm) for details.  We consider here only the special case
that~$(\UU_t)$ is an i.i.d.\ sequence, i.e.,\ that $q(v, w) \equiv q(w)$.  A transition
rule~$\phi$ for the chain~$(\XX_t)$ is created as follows:\ for vertex~$u$ and
arborescence~$x$ with root~$r$, $\phi(x, u)$ is the arborescence obtained from~$x$
by adding an arc from~$r$ to~$u$ and deleting the unique arc in~$x$ whose tail is~$u$.
Then it can be shown that if $\DD_s$ is taken to be the set of all vertices appearing
at least once in $(\UU_1, \ldots, \UU_s)$ and
$\Delta := \{\Uc\}$, then~$\vecDD$ is a detection process.
}
\end{example}

\subsection{Bounding processes}
\label{bounding}

\hspace{\parindent}
We obtain a natural and useful example of a detection process~$\vecDD$ when
(a)~$\vecDD$ is constructed solely from~$\vecUU$,
(b)~the corresponding state space~$\Dc$ is some
collection of subsets of~$\Xc$, with
$$
\Delta := \{ \{z\}:\ z \in \Xc\},
$$
and
$$
\mbox{(c)\ \ \ }\DD_s \supseteq \{\YY_s(x):\ x \in \Xc\}.
$$
The concept is simple:\ in this case, each set~$\DD_s$ is just a ``conservative
estimate'' (i.e., a superset) of the corresponding set $\{\YY_s(x): x \in \Xc\}$
of trajectory values; thus if $\DD_s = \{z\}$, then the trajectories~$\vecYY(x)$
are coalesced to state~$z$ at time~$s$ and remain coalesced thereafter.  We
follow the natural impulse to call such a set-valued detection process
a {\em bounding process\/}.  Such bounding processes arise naturally in the contexts of
monotone (and anti-monotone) transition rules, and have been used by many authors:\
see Section~\ref{mono}.  Other examples of bounding processes can be found in
H\"{a}ggstr\"{o}m and Nelander~\cite{HN2} (for CFTP) and in works
of Huber:\ see~\cite{Hefficient} and~\cite{Hexact} in connection with CFTP
and~\cite{Hinterruptible} in connection with our algorithm.

Of course, nothing is gained, in comparison to tracking all the trajectories, by
the use of a bounding process unless the states of~$\Dc$ have more concise
representations than those of generic subsets of~$\Xc$; after all,
we could always choose $\Dc = 2^{\Xc}$ and $\DD_s =
\{\YY_s(x): x \in \Xc\}$.  One rather general, and frequently applied,
setting where compact representations are possible, discussed in Section~\ref{mono},
is that of a realizably monotone chain on a partially ordered set (poset)~$\Xc$.
See the references listed in Remark~\ref{anti}(a) for examples of the use
of bounding processes for such chains.

See also Remark~\ref{anti}(b) for a brief discussion of the closely related
concept of ``dominating process.''

\begin{remark}
\label{detvsbd}
{\em
With our algorithm, one can take advantage of detection processes that are more general than
bounding processes.   Consider Example~\ref{MTFandtrees}(a) as an illustrative case.
Suppose we start the first (i.e.,\ time-reversed) phase of our algorithm in the
identity permutation (as is done in~\cite{MTFalgo}).  When we run the forward-time
phase of the algorithm, we only need to keep track of the {\em unordered\/} set of
records requested; this is what the detection process~$\vecDD$ of
Example~\ref{MTFandtrees}(a) does.  If the event~$\EE$ of Section~\ref{detectproc}
occurs, that is, if all, or all but one, of the records have been requested at least
once by time~$t$, then we have detected coalescence.  Then, without further work, we
know the state to which the trajectories have coalesced:\ namely, our initial identity
permutation.  To maintain a bounding process for this same example, we would have to do
more work:\ we would have to keep track of the {\em ordered\/} set of records requested,
ordered by most recent request.  (Note that a bounding process {\em could\/} be constructed
by combining information from the dominating process and the trajectory generated in the
algorithm's initial backward phase.  But it is not useful to do so.)
}
\end{remark}

\sect{A mathematically rigorous framework}
\label{frame}

\subsection{The formal framework}
\label{formal}

\hspace\parindent
We now formally set up a general framework for rejection sampling using auxiliary
randomness, paying careful attention to technical details.  We will apply this framework
not only to provide a rigorous treatment of Algorithm~\ref{nontech} (in Section~\ref{app1}),
but also (independently) to handle a variant (Algorithm~\ref{SMnontech}) that applies
to a stochastically monotone kernel~$K$.

We need to set up two probability spaces.  In our main application to
Algorithm~\ref{nontech}, the first space $(\Omega, \Ac, P)$---designated in
ordinary typeface---will be useful for theoretical considerations and for the
computation of certain conditional probability distributions.  The second
space $(\OO, \AAA, \PP)$---designated in boldface type---will be the probability
space actually simulated when the algorithm is run.
All random variables defined on the first
space (respectively, second space) will also be designated in ordinary typeface (resp.,\
boldface type). We have chosen this notational system to aid the reader:\ Corresponding
variables, such as~$X_0$ and~$\XX_0$, will play analogous roles in the two spaces.

Recall that for a measurable
space $(\Zc, \Fc)$ and a (not necessarily measurable) subset $E \subseteq \Zc$, the
{\em trace $\sigma$-field\/} (on the sample space~$E$) is the $\sigma$-field $\{E
\cap F: F \in \Fc\}$.  In setting up both probability spaces, we assume:
\begin{enumerate}
\item[(i)] $(\Xc, \Bc)$ and $(\Xc', \Bc')$ are measurable spaces,  $B' \subseteq
\Xc'$, and $\Bc''$ is the trace $\sigma$-field for~$B'$
\item[(ii)] $\pi$ is a probability measure on~$\Bc$.
\end{enumerate}
For our first probability space, we also assume that
\begin{enumerate}
\item[(iii)] $(\Omega, \Ac, P)$ is a probability space on which are defined
mappings $X: \Omega \to \Xc$ and $X': \Omega \to \Xc'$ and $Y: \Omega \to \Xc'$
and a set $C \subseteq \Omega$
\item[(iv)] $X$ is measurable from $\Ac$ to $\Bc$, with $\Lc(X) = PX^{-1}
= \pi$
\end{enumerate}
and that the following two conditions hold, for every $B'' \in \Bc''$:
\begin{enumerate}
\item[(v)] $\{X' \in B''\} \cap C \in \Ac$ and $\{Y \in B''\} \cap C \in \Ac$
\item[(vi)] $P(\{X' \in B''\} \cap C\,|\,X \in B) = P(\{Y \in B''\} \cap C)$
for all $B \in \Bc$ with $\pi(B) > 0$.
\end{enumerate}
Notice that~(v) and~(vi) are implied by
\begin{enumerate}
\item[(v$'$)] $\{X' \in B''\} \cap C \in \Ac$,\ \ $Y = X'$ on the set~$C$, and~$X$
and the event $\{Y \in B''\} \cap C$ are independent.
\end{enumerate}

\begin{proposition}
\label{basicprop}
Under assumptions {\rm (i)--(vi)},
\begin{equation}
\label{basicfact}
P(\{X' \in B''\} \cap C \cap \{X \in B\}) = P(\{X' \in B''\} \cap
C)\,\pi(B)
\end{equation}
for every $B \in \Bc$ and every $B'' \in \Bc''$.
\end{proposition}

\proof
We may assume $\pi(B) > 0$, in which case~(vi) implies
\begin{eqnarray}
\label{facty}
P(\{X' \in B''\} \cap C \cap \{X \in B\})
 &=& P(\{Y \in B''\} \cap C)\,\pi(B),\\
\label{sub}
P(\{X' \in B''\} \cap C)
 &=& P(\{Y \in B''\} \cap C)
\end{eqnarray}
for general~$B \in \Bc$ and for~$B = \Xc$, respectively.  Substituting~(\ref{sub})
into~(\ref{facty}) gives~(\ref{basicfact}).~\qed

\bigskip
Two corollaries follow immediately:

\begin{corollary}
\label{monocor}
Under assumptions {\rm (i)--(vi)}, if $P(\{X' \in B'\} \cap C) > 0$, then
$$
\Lc(X\,|\,\{X' \in B'\} \cap C) = \pi.
$$
\end{corollary}

\begin{corollary}
\label{generalcor}
Suppose assumptions {\rm (i)--(vi)} hold for $B'  = \Xc'$; in
particular, {\rm (v)} then implies that $C \in \Ac$.  Suppose also that $X'$ is
measurable from $\Ac$ to $\Bc'$.  Assume $P(C) > 0$.  Then
$$
\mbox{\rm $X'$ and~$X$ are conditionally independent given~$C$, and $\Lc(X | C) = \pi$.}
$$
\end{corollary}
\bigskip

Now we set up the second probability space.  Specifically, consider the following
assumptions: that
\begin{enumerate}
\item[(vii)] $(\OO, \AAA, \PP)$ is a probability space on which are defined
mappings $\XX: \OO \to \Xc$ and $\XXp: \OO \to \Xc'$ and a set $\CC \subseteq \OO$
\item[(viii)] $\XX$ is measurable from $\AAA$ to $\Bc$,
\end{enumerate}
that, for every $B'' \in \Bc''$, we have
\begin{enumerate}
\item[(ix)] $\{\XXp \in B''\} \cap \CC \in \AAA$,
\end{enumerate}
and that we have the following basic connection between the two spaces:
\begin{enumerate}
\item[(x)] The measure
$$
\PP(\{\XXp \in dx'\} \cap \CC \cap \{\XX \in dx\})
$$
on the product space $(B' \times \Xc, \Bc'' \otimes \Bc)$ has a density $D(x, x')
\equiv D(x')$ ($x' \in B'$) that doesn't depend on $x \in \Xc$ with respect to the
measure 
$$
P(\{X' \in dx'\} \cap \ C \cap \{X \in dx\}).
$$ 
\end{enumerate}
Notice that~(x) is implied by conditions (x$'$)--(x$'''$),
wherein $\Lc(\XXp) = \PP (\XXp)^{-1}$:
\begin{enumerate}
\item[(x$'$)] $B' = \Xc'$, and $X'$ and~$\XXp$ are measurable (from~$\Ac$
and~$\AAA$, respectively, to~$\Bc'$)
\item[(x$''$)] $\Lc(\XXp) \ll \Lc(X')$, with Radon--Nikodym derivative~$D$
\item[(x$'''$)] There exists a conditional subprobability distribution
$$
P(C \cap \{X \in dx\}\,|\,X' = x'),\ \ x' \in \Xc',
$$
which also serves as conditional subprobability distribution
$$
\PP(\CC \cap \{\XX \in dx\}\,|\,\XXp = x'),\ \ x' \in \Xc'.
$$
\end{enumerate}

It is now key that the results of Proposition~\ref{basicprop} and
Corollaries~\ref{monocor} and~\ref{generalcor} carry over to our second space: 

\begin{proposition}
\label{bbasicprop}
Under assumptions {\rm (i)--(x)},
\begin{equation}
\label{bbasicfact}
\PP(\{\XXp \in B''\} \cap \CC \cap \{\XX \in B\}) = \PP(\{\XXp \in B''\} \cap
\CC)\,\pi(B)
\end{equation}
for every $B \in \Bc$ and every $B'' \in \Bc''$.
\end{proposition}

\proof
Let $D$ be the Radon--Nikodym derivative guaranteed by assumption~(x).  Then, using
Proposition~\ref{basicprop},
\begin{eqnarray*}
\PP(\{\XXp \in B''\} \cap \CC \cap \{\XX \in B\})
 &=& \int_{B''}\,D(x')\,P(\{X' \in dx'\} \cap C \cap \{X \in B\}) \\
 &=& \int_{B''}\,D(x')\,P(\{X' \in dx'\} \cap C)\,\pi(B).
\end{eqnarray*}
As usual, setting $B = \Xc$ and substituting, we obtain~(\ref{bbasicfact}).~\qed

\begin{corollary}
\label{mmonocor}
Under assumptions {\rm (i)--(x)}, if $\PP(\{\XXp \in B'\} \cap \CC) > 0$, then
$$
\Lc(\XX\,|\,\{\XXp \in B'\} \cap \CC) = \pi.
$$
\end{corollary}

\begin{corollary}
\label{ggeneralcor}
Suppose assumptions {\rm (i)--(x)} hold for $B'  = \Xc'$.  Suppose also that $\XXp$
is measurable from $\AAA$ to $\Bc'$.  Assume $\PP(\CC) > 0$.  Then
$$
\mbox{\rm $\XXp$ and~$\XX$ are conditionally independent given $\CC$,
and $\Lc(\XX|\CC) = \pi$.}
$$
\end{corollary}

\subsection{Details for Algorithm~2.1}
\label{app1}

\label{rigor}

\hspace\parindent
The goal of this subsection is to describe when and how Algorithm~\ref{nontech} can
be applied legitimately.
\medskip

{\em The space $(\Xc, \Bc)$:\/}\ Recall that a {\em Polish space\/} is a complete
separable metric space.  For convenience, we shall assume that the measurable
state space $(\Xc, \Bc)$ of interest (on which the probability measure~$\pi$ of
interest is defined) is isomorphic to a Borel subset of a Polish space (with its
trace Borel $\sigma$-field).  This assumption will at once guarantee the existence
of such objects as conditional distributions that would otherwise require
individual arguments or assumptions.  We call such a space a {\em standard Borel
space\/}.  Our assumption should cover most cases of applied interest.
\medskip

{\em The kernel~$K$ and its time-reversal~$\Kt$:\/}\ Let~$K: \Xc \times \Bc \to
[0,1]$ be a Markov transition kernel on~$\Xc$; that is, we suppose
that~$K(x,\cdot)$ is a probability measure on~$\Bc$ for each $x \in \Xc$ and that
$K(\cdot, B)$ is a $\Bc$-measurable function for each~$B \in \Bc$.  The kernel is
chosen (by the user) so that~$\pi$ is a stationary distribution, i.e.,\ so that
$$
\int_{\Xc} \pi(dx) K(x, dy) = \pi(dy)\mbox{\ \ \ on $\Xc$.}
$$
Since~$(\Xc, \Bc)$ is standard Borel, there exists ``a conditional distribution the
other way around''---more precisely, a Markov kernel~$\Kt$ on~$\Xc$ satisfying
\begin{equation}
\label{reversal}
\pi(dx) K(x, dy) = \pi(dy) \Kt(y, dx)\mbox{\ \ \ on $\Xc \times \Xc$.}
\end{equation}
Given~$\pi$ and~$K$, the kernel~$\Kt(y, dx)$ is $\pi(dy)$-almost surely uniquely
defined.  We choose and fix such a~$\Kt$.
\medskip

{\em The transition rule~$\phi$:\/}\ It can be shown that there exists a transition
rule which can be used to drive the construction of the Markov chain of interest. 
More precisely, our assumption that~$(\Xc, \Bc)$ is standard Borel implies that
there exists a standard Borel space~$(\Uc, \Fc)$, a product-measurable function
$\phi: \Xc \times \Uc \to \Xc$, and a probability measure~$\mu$ on~$\Fc$, such that
\begin{equation}
\label{kq}
K(x, B) = \mu\{u: \phi(x, u) \in B\},\ \ \ x \in \Xc,\ \ B \in \Bc.
\end{equation}
Such~$\phi$ (with accompanying~$\mu$) is sometimes called a {\em transition
rule\/}.  We choose and fix such a $(\phi, \mu)$.

\begin{remark}
\label{phiremark}
{\em
(a)~Conversely, if~$\phi$ has the stated properties and~$K$ is {\em defined\/}
by~(\ref{kq}), then~$K$ is a Markov kernel.

(b)~A transition rule~$\phi$ can always be found that uses $(\Uc, \Fc, \mu) =
([0, 1], \mbox{Borels},$ $\mbox{uniform distribution})$.  The proof of existence
(cf.\ Theorem~1.1 in Kifer~\cite{Kifer} and Remark~(iv) at the end of Section~5.2
in Diaconis and Freedman~\cite{DF}) makes use of inverse probability transforms and
certain standard reduction arguments.  In the special case that~$(\Xc, \Bc) =
([0,1], \mbox{Borels})$, we can in fact use
$$
\phi(x, u) \equiv G(x, u)
$$
where~$G(x, \cdot)$ is the usual inverse probability transform corresponding to
the distribution function $u \mapsto K(x, [0,u])$.

(c)~If~$(\Xc, \Bc)$ is any discrete space (i.e.,\ if~$\Xc$ is countable and~$\Bc$
is the total $\sigma$-field), a very simple alternative choice is the following
``independent-transitions'' transition rule.  Let~$\Uc = \Xc^{\Xc}$ (with~$\Fc$ the
product $\sigma$-algebra), let $\mu$ be product measure with $x$th
marginal $K(x, \cdot)$ ($x \in \Xc$), and let~$\phi$ be the evaluation function
$$
\phi(x, u) := u(x).
$$

(d) Many interesting examples of transition rules can be found in the literature,
including Diaconis and Freedman~\cite{DF} and the references cited in
Section~\ref{intro}.

(e) Usually there is a wealth of choices of transition rule, and the art is to
find one giving rapid and easily detected coalescence.  Without going into details
at this point, we remark that the transition rule in~(c) usually performs quite
badly, while transition rules having a certain monotonicity property will perform
well under monotonicity assumptions on~$K$.
}
\end{remark}
\medskip

{\em The Markov chain and the first probability space:\/}\ From our previous
comments it is now easy to see that there exists a standard Borel space~$(\Uc,
\Fc)$, a transition rule $(\phi, \mu)$, and a probability space~$(\Omega,
\Ac, P)$ on which are defined independent random variables $X_0, U_1, U_2, \ldots,
U_t$ with $X_0 \sim \pi$ and each $U_s \sim \mu$.  Now inductively define
\begin{equation}
\label{drive}
X_s := \phi(X_{s - 1}, U_s),\ \ \ 1 \leq s \leq t.
\end{equation}
Then $\vecX := (X_0, \ldots, X_t)$ is easily seen to be a stationary Markov
chain with kernel~$K$, in the sense that 
\begin{equation}
\label{ftraj}
\ \ \ \ \ \ \ P(X_0 \in dx_0, \ldots, X_t \in dx_t) = \pi(dx_0) K(x_0,
dx_1)\!\cdots\!K(x_{t-1}, dx_t)\mbox{\ \ \ on~$\Xc^{t+1}$.}
\end{equation}
In fact, for each $x \in \Xc$ we obtain a chain with kernel~$K$ started from~$x$
by defining $Y_0(x) := x$ and, inductively,
$$
Y_s(x) := \phi(Y_{s-1}(x), U_s).
$$
Let $\vecY(x) := (Y_0(x), \ldots, Y_t(x))$.  In this notation we have
$\vecY(X_0) = \vecX$.  Recalling the notational note at the end of
Section~\ref{intro},
let
\begin{equation}
\label{Cdef}
C := \{ Y_t(x) \mbox{\ does not depend on~$x$}\}
\end{equation}
denote the set of sample points~$\omega$ for which the trajectories~$\vecY(x)$
have all coalesced by time~$t$.  We assume that~$C$ belongs to the
$\sigma$-field $\sigma \langle \vecU \rangle$ generated by $\vecU := (U_1, \ldots, U_t)$.

\begin{remark}
{\em
For this remark, suppose that~$\Xc$ {\em is\/} a Borel subset of a Polish space
(and hence a separable metric space in its own right).
We will prove that continuity of the transition rule $\phi(x, u)$ in~$x \in \Xc$
for each $u \in \Uc$ is sufficient for $C \in \sigma \langle \vecU \rangle$,
and we note that this is automatic if~$\Xc$ is discrete.

As guaranteed by separability of~$\Xc$, let~$D$ be a countable dense subset of~$\Xc$.
Given $z \in \Xc$ and $\epsilon \geq 0$, let $B_z(\epsilon)$ denote the closed ball of
radius~$\epsilon$ centered at~$z$.  

Suppose that~$\phi(\cdot, u)$ is continuous for each $u \in \Uc$.  Define the
iterates $\phi^s:~\Xc~\times \Uc^s \to \Xc$ ($s = 1, 2, \ldots, t$)
inductively by $\phi^1 := \phi$ and
$$
\phi^s(x; u_1, \ldots, u_s) := \phi(\phi^{s - 1}(x; u_1, \ldots, u_{s - 1}); u_s).
$$
Note that~$\phi^t$ is, like~$\phi$, continuous in its first argument,
and that $Y_t(x) = \phi^t(x; \vec{U})$.
Using the separability of~$\Xc$, it is not hard to show that
$$
C = \cap^{\infty}_{n = 1} \cup_{z \in D} \cap_{x \in D}
\{\phi^t(x; \vec{U}) \in B_z(1 / n)\},
$$
from which the desired measurability of~$C$ is evident.
}
\end{remark}

Now observe that conditions (i)--(iv) and (v$'$) in Section~\ref{formal} are
satisfied by fixing $x^*_0 \in \Xc$ arbitrarily and taking
\begin{eqnarray}
\label{defs}
(\Xc', \Bc') &=& (\Xc, \Bc),\ \ \ B' = \Xc',\ \ \ \Bc'' = \Bc, \\
\nonumber
X &=& X_0,\ \ \ X' = X_t,\ \ \ Y = Y_t(x^*_0),\ \ \ \mbox{$C$ as at~(\ref{Cdef}).}
\end{eqnarray}
Note that the independence in~(v$'$) follows from the fact that~$X_0$
and~$\vecU$ have been chosen to be independent.
\medskip

{\em The second probability space and the algorithm:\/}\ The key to setting up
the second probability space is to satisfy assumption~(x$'''$) in
Section~\ref{formal}.  In calculating the first-space conditional distribution
mentioned there, we will make use of the auxiliary randomness provided by $X_1,
\ldots, X_{t - 1}$ and $\vecU$ and compute in stages.  First observe
from~(\ref{ftraj}) and repeated use of~(\ref{reversal}) that
$$
P(X_0 \in dx_0, \ldots, X_{t - 1} \in dx_{t - 1}\,|\,X_t = x_t) = \Kt(x_t, dx_{t -
1}) \cdots \Kt(x_1, dx_0)
$$
serves as a conditional distribution for $(X_0, \ldots, X_{t - 1})$ given $X_t =
x_t$.  Next, we will discuss in Section~\ref{imp} how to compute
$\Lc(\vecU\,|\,\vecX = \vecx)$.  Finally, from our assumption that $C \in
\sigma(\vecU)$, it follows that we can write the indicator $Z := I_C$ as $Z =
\Gamma(\vecU)$ for some product-measurable $\Gamma: \Uc^t \to \{0, 1\}$, and one
can check the intuitively obvious assertion that unit mass at $\Gamma(\vecu)$
serves as a conditional distribution for~$Z$ given $(\vecX, \vecU) = (\vecx,
\vecu)$.  We get the conditional distribution in~(x$'''$) by chaining together the
conditional distributions we have computed and integrating out the auxiliary
variables, in the obvious and standard fashion.
 
Moreover, our discussion has indicated how to set up and {\em simulate\/} the
second space.  To satisfy assumption~(x$''$) in Section~\ref{formal}, we assume
that the law of $\XX_t$ 
chosen by the user
is absolutely continuous with respect to $\pi$; of course we do {\em not\/}
assume that the user can compute the Radon--Nikodym derivative~$D$.  [For example,
in the common situation that~$(\Xc, \Bc)$ is discrete and $\pi({x}) > 0$ for every
$x
\in
\Xc$, the value of~$\XX_t$ can be chosen deterministically and arbitrarily.] 
Having chosen $\XX_t = x_t$, the user draws an observation $\XX_{t - 1} = x_{t -
1}$ from
$\Kt(x_t,
\cdot)$, then an observation $\XX_{t - 2} = x_{t - 2}$ from $\Kt(x_{t - 1},
\cdot)$, etc.\ \ Next, having chosen $\vecXX = \vecx$ [i.e., $(\XX_0, \ldots,
\XX_t) = (x_0, \ldots, x_t)$], the user draws an observation $\vecUU = \vecu$
from $\Lc(\vecU\,|\,\vecX = \vecx)$.  Finally, the user sets $\ZZ = \Gamma(\vecu)$
and declares that~$\CC$, or {\em coalescence\/}, has occurred if and only if $\ZZ =
1$.  With the definitions~(\ref{defs}), $\CC$ as above, and
$$
\XX = \XX_0,\ \ \ \XXp = \XX_t,
$$
assumptions (i)--(x) are routinely verified.  According to
Corollary~\ref{ggeneralcor}, if \linebreak $\PP(\CC) > 0$,
then $\Lc(\XX_0|\CC) = \pi$.  It follows that the conditional distribution
of output from Algorithm~\ref{nontech} given that it ultimately succeeds (perhaps
only after many iterations of the basic routine) is~$\pi$, as desired.

\begin{remark}
\label{genremark}
{\em (a)~If $\PP(\CC) > 0$ for suitably large~$t$, then ultimate success is (a.s.)\
guaranteed if the successive choices of~$t$ become large.  A necessary condition
for ultimate positivity of~$\PP(\CC)$ is uniform ergodicity of~$K$.  This condition
is also sufficient, in the (rather weak) sense that if~$K$ is uniformly ergodic,
then there exists a finite integer~$m$ and a transition rule~$\phi_m$ for
the $m$-step kernel~$K^m$ such that Algorithm~\ref{nontech}, applied
using~$\phi_m$, has $\PP(\CC) > 0$ when~$t$ is chosen sufficiently large.
Compare the analogous Theorem~4.2 for CFTP in Foss and Tweedie~\cite{FT}.

A similar remark applies to Algorithm~\ref{altalg}.

(b)~Just as discussed in Fill~\cite{Fill} (see especially the end of Section~7
there), the algorithm (including its repetition of the basic routine) we have
described is interruptible; that is, its running time (as measured by number of
Markov chain steps) and output are independent random variables, conditionally
given that the algorithm eventually terminates.

(c)~If the user chooses the value of~$\XX_t$ ($= z$, say) deterministically, then
all that can be said in general is that the algorithm works properly for
$\pi$-a.e.\ such choice.  In this case, let the notation $\PP_z(\CC)$ reflect the
dependence of~$\PP(\CC)$ on the initial state~$z$.  Then clearly
$$
\int\,\PP_z(\CC)\,\pi(dz) = P(C),
$$ 
which is the unconditional probability of coalescence in our first probability
space and therefore equal to the probability that CFTP terminates over an interval
of width~$t$.  This provides a first link between CFTP and
Algorithm~\ref{nontech}.  Very roughly recast, the distribution of running time
for CFTP is the stationary mixture, over initial states, of the distributions of
running time for Algorithm~\ref{nontech}.  For further elaboration of the
connection between the two algorithms, see Section~\ref{conn}.
}
\end{remark}

\subsection{Imputation}
\label{imp}

\hspace\parindent
In order to be able to run Algorithm~\ref{nontech}, the user needs to be able to
impute~$\vecU$ from~$\vecX$, i.e.,\ to draw from $\Lc(\vecU\,|\,\vecX = \vecx)$. 
In this subsection we explain how to do this.

We proceed heuristically at first:
\begin{eqnarray*}
\lefteqn{P(\vecU \in d\vecu\,|\,\vecX = \vecx)} \\
 & =& P(\vecU \in d\vecu\,|\,X_0 = x_0,\,\phi(x_0, U_1) = x_1,\,\ldots,\,\phi(x_{t
        - 1}, U_t) = x_t)\mbox{\ \ \ by~(\ref{drive})} \\
 & =& P(\vecU \in d\vecu\,|\,\phi(x_0, U_1) = x_1,\,\ldots,\,\phi(x_{t - 1}, U_t) =
        x_t)\mbox{\ \ \ by indep.\ of~$X_0$ and~$\vecU$} \\
 & =& P(U_1 \in du_1\,|\,\phi(x_0, U_1) = x_1) \times \cdots \times P(U_t \in
        du_t\,|\,\phi(x_{t - 1}, U_t) = x_t) \\
 &{}& \mbox{\ \ \ \ \ \ \ by independence of $U_1, \ldots, U_t$} \\
 & =& P(U_1 \in du_1\,|\,\phi(x_0, U_1) = x_1) \times \cdots \times P(U_1 \in
        du_t\,|\,\phi(x_{t - 1}, U_1) = x_t) \\
 &{}& \mbox{\ \ \ \ \ \ \ since $U_1, \ldots, U_t$ are identically distributed} \\
 & =& P(U_1 \in du_1\,|\,X_0 = x_0,\,X_1 = x_1) \times \cdots \times P(U_1 \in
        du_t\,|\,X_0 = x_{t - 1},\,X_1 = x_t),
\end{eqnarray*}
where the last equality is justified in the same fashion as for the first two.

In fact, the result of this heuristic calculation is rigorously correct; its
proof is an elementary but not-entirely-trivial exercise in the use of conditional
probability distributions.  The existence (and a.s.\ uniqueness) of a conditional
distribution $\Lc(U_1\,|X_0 = \cdot,\,X_1 = \cdot)$ is guaranteed by the fact
that~$(\Uc, \Fc)$ is standard Borel; moreover,
\begin{lemma}
\label{condlemma}
The $t$-fold product of the measures
$$
P(U_1 \in du_1\,|\,X_0 = x_0,\,X_1 = x_1),\ \ldots,\  P(U_1 \in du_t\,|\,X_0 = x_{t
  - 1},\,X_1 = x_t)
$$
serves as a conditional probability distribution $P(\vecU \in d\vecu\,|\,\vecX =
\vecx)$.
\end{lemma}

In setting up the second probability space, therefore, the user, having chosen
$\vecXX = \vecx$, draws an observation $\vecUU = \vecu$ by drawing $\UU_1, \ldots,
\UU_t$ independently, with $\UU_s$ chosen according to the distribution
$\Lc(U_1\,|\,X_0 = x_{s - 1},\,X_1 = x_s)$.

\begin{remark}
\label{subtle}
{\em (a)~There are subtleties involved in the rigorous proof of
Lemma~\ref{condlemma}.  In particular, there is no justification apparent to us,
in general, that the conditional distributions $\Lc(U_1\,|\,\phi(x_0, U_1) =
\cdot)$, one for each fixed $x_0 \in \Xc$, can be chosen in such a way that
$P(U_1 \in F\,|\,\phi(x_0, U_1) = x_1)$ is jointly measurable in~$(x_0, x_1)$
for each $F \in \Fc$.  Nevertheless, if we {\em assume\/} such measurability, then
one can show rigorously that
$$
P(U_1 \in du_1\,|\,\phi(x_0, U_1) = x_1) \times \cdots \times P(U_1 \in
du_t\,|\,\phi(x_{t - 1}, U_1) = x_t)
$$
serves as $P(\vecU \in d\vecu\,|\,\vecX = \vecx)$.

(b)~If~$(\Xc, \Bc)$ is discrete, then of course the measurability in~(a) is
automatic.  Suppose we use the ``independent-transitions'' rule~$\phi$ discussed
in Remark~\ref{phiremark}(c).  Then the measure~$\mu$, but with the $x_0$th
marginal replaced by~$\delta_{x_1}$, serves as $\Lc(U_1\,|\,\phi(x_0, U_1) = x_1) =
\Lc(U_1\,|\,U_1(x_{0}) = x_1)$ and therefore as $\Lc(U_1\,|\,X_0 = x_0,\,X_1
= x_1)$.  Informally stated, having chosen $\XX_s = x_s$ and $\XX_{s - 1} = x_{s -
1}$, the user imputes the forward-trajectory transitions from time~$s - 1$ to
time~$s$ in Algorithm~\ref{nontech} by declaring that the transition from
state~$x_{s - 1}$ is to state~$x_s$ and that the transitions from other states are
chosen independently according to their usual non-$\vecX$-conditioned
distributions.

(c)~As another example, suppose that $\Xc = [0, 1]$ and we use the inverse
probability transform transition rule discussed in Remark~\ref{phiremark}(b). 
Suppose also that each distribution function $F(x_0, \cdot) = K(x_0, [0, \cdot])$
is strictly increasing and onto~$[0, 1]$ and that $F(x_0, x_1)$ is jointly
Borel-measurable in~$x_0$ and~$x_1$.  Then~$\delta_{F(x_0, x_1)}$ serves as
$\Lc(U_1\,|\,X_0 = x_0,\,X_1 = x_1)$.  Informally stated, a generated pair
$(\XX_s, \XX_{s - 1}) = (x_s, x_{s-1})$ completely determines the value $F(x_{s -
1}, x_s)$ for~$\UU_s$.
}
\end{remark}

\sect{Monotonicity}
\label{app3}

\hspace{\parindent}
Throughout Section~\ref{app3} we suppose that
$(\Xc, \Bc)$ is a Polish space with (Borel $\sigma$-field
and) closed partial order $\leq$; the meaning of {\em closed\/} here is that
$\{(x, y) \in \Xc \times \Xc: x \leq y\}$ is assumed to be closed in the
product topology.  [For example, closedness is automatic for any partial order if
$(\Xc, \Bc)$ is discrete.]
We also assume that there exist (necessarily unique) elements~$\zh$ and~$\oh$
in~$\Xc$ (called {\em bottom element\/} and {\em top element\/}, respectively) such that
$\zh \leq x \leq \oh$ for all $x \in \Xc$.

We require the notion of stochastic monotonicity, according to
the following definitions (which extend Definition~4.1 of Fill~\cite{Fill}
to our more general setting).

\begin{definition}
\label{stochdef}
{\em
(a)~A subset~$B$ of~$\Xc$ is called a {\em down-set\/} or {\em order ideal\/}
if, whenever $x \in B$ and $y \leq x$, we have $y \in B$.

(b)~Given two probability measures~$\nu_1$ and~$\nu_2$ on~$\Bc$, we
say that {\em $\nu_1 \leq \nu_2$ stochastically\/}, and write
$\nu_1 \preceq \nu_2$, if $\nu_1(B)
\geq \nu_2(B)$ for every closed down-set~$B$.

(c)~A kernel~$K$ is said to be {\em stochastically monotone\/} (SM) if $K(x,
\cdot) \preceq K(y, \cdot)$ whenever $x \leq y$.
}
\end{definition}

The main goal of this section is to describe an analogue of
Algorithm~\ref{nontech} which applies when only stochastic monotonicity of~$K$ is
assumed.  Here is a rough formulation; the details will be discussed in
Section~\ref{SMrigor}.

\begin{algorithm}
\label{SMnontech}
{\em
Consider a stochastically monotone kernel $K$ on a partially ordered set with
bottom element~$\zh$ and top element~$\oh$.
Choose and fix a positive integer~$t$, set $\XX_t = \zh$, and perform the
following routine.  Run the time-reversed chain~$\Kt$ for~$t$ steps, obtaining
$\XX_t, \XX_{t - 1}, \ldots, \XX_0$ in succession.  Then, reversing the direction
of time, generate a chain, say $\vecYY = (\YY_0, \ldots, \YY_t)$, with $\YY_0 =
\oh$ and kernel~$K$; this trajectory is to be coupled {\em ex post facto\/} with
$\vecXX = (\XX_0, \ldots, \XX_t)$, which is regarded as a trajectory from~$K$. 
Finally, we check whether $\YY_t = \zh$.  If so, the value~$\XX_0$ is accepted as
an observation from~$\pi$; if not, we repeat (as for Algorithm~\ref{nontech},
but always with~$\XX_t = \zh$).
}
\end{algorithm}

\begin{remark}
{\em
When~$\Xc$ is finite, Algorithm~\ref{SMnontech} reduces to the algorithm of
Section~7.2 in Fill~\cite{Fill}.
}
\end{remark}

\subsection{Realizable monotonicity}
\label{mono}

\hspace{\parindent}
It is easy to see from~(\ref{kq}) that if there exists
a monotone transition rule, i.e.,\ a transition rule~$\phi$ 
with the property that
\begin{equation}
\label{mondef}
\phi(x, u) \leq \phi(y, u)\mbox{\ \ for every~$u$}\quad {\rm whenever} \quad x \leq y,
\end{equation}
then~$K$ is stochastically monotone according to Definition~\ref{stochdef}(c).
We call this stronger property {\em realizable monotonicity\/}.
It is a common misbelief that, conversely, stochastic
monotonicity for~$K$ implies realizable monotonicity.
This myth is annihilated by Fill and Machida~\cite{SMRM} and Machida~\cite{Machida},
even for the case of a finite poset~$\Xc$,
for which it is shown that
every stochastically monotone kernel is
realizably monotone (i.e.,\ admits a monotone transition rule) if and only if the
cover graph of~$\Xc$ (i.e.,\ its Hasse diagram regarded as an undirected graph) is
acyclic.

Nevertheless, realizable monotonicity can be used to motivate and explain
Algorithm~\ref{SMnontech}.  Indeed,
suppose for the remainder of this subsection that~$K$ admits a monotone
transition rule~$\phi$
as in~(\ref{mondef}).  Then we proceed to build a bounding process as in
Section~\ref{bounding} and show how Algorithm~\ref{nontech} can be applied
efficiently.
One immediately verifies by induction that
\begin{equation}
\label{squeeze}
\YY_s(\zh) \leq \YY_s(x) \leq \YY_s(\oh)\mbox{\ \ \ 
for all $0 \leq s \leq t$ and all $x \in \Xc$.}
\end{equation}
Thus $\DD_s := [\YY_s(\zh), \YY_s(\oh)] = \{y \in \Xc:\ \YY_s(\zh) \leq y \leq \YY_s(\oh)\}$
gives a bounding process, and the pair $(\YY_s(\zh), \YY_s(\oh))$ is a quite concise
representation of~$\DD_s$.  In plain language, since
monotonicity is preserved, when the chains~$\vecYY(\zh)$ and~$\vecYY(\oh)$ have
coalesced, so must have every~$\vecYY(x)$.

But note also that, if we choose the initial state~$\XX_t$ to be~$\zh$,
then~$\{\YY_t(\oh) = \zh\}$ is the coalescence event~$\CC$.
Algorithmically, it
follows that if $\XX_t = \zh$ is a legitimate starting point for
Algorithm~\ref{nontech} [as discussed near the end of Section~\ref{rigor}, it is
sufficient for this that $\pi(\{\zh\}) > 0$], and if $\PP(\CC) > 0$,
then $\Lc(\XX_0|\CC) = \pi$.  Informally put, we need only track the
single upper-bound trajectory~$\vecYY(\oh)$ in the forward phase; if $\YY_t(\oh) = \zh$,
then the routine (correctly) accepts~$\XX_0$ as an observation from~$\pi$.

Notwithstanding the fundamental distinction
between the two notions of monotonicity, the routine of Algorithm~\ref{SMnontech}
shares with its cousin from the realizably monotone case the feature that, in the
forward phase, the user need only check whether a single $K$-trajectory started
at~$\oh$ ends at~$\zh$.

\begin{remark}
\label{anti}
{\em (a)~Lower and upper bounding processes can also be constructed when
Algorithm~\ref{nontech} is applied with a so-called ``anti-monotone'' transition
rule; we omit the details.  See H\"{a}ggstr\"{o}m and Nelander~\cite{HN},
Huber~\cite{Hexact}, Kendall~\cite{Kendall}, M{\o}ller~\cite{Moller}, M{\o}ller
and Schladitz~\cite{MS}, and Th\"{o}nnes~\cite{Thonnes} for further discussion in
various specialized settings.  There are at least two neat tricks associated with
anti-monotone rules.  The first is that, by altering the natural partial order
on~$\Xc$, such rules can be regarded, in certain bipartite-type settings, as
monotone rules, in which case the analysis of Section~\ref{performance} 
(with~$\zh$ and~$\oh$ taken in the {\em altered\/} ordering, of course) is
available:\ consult Section~3 of~\cite{HN}, the paper~\cite{MS}, and Definition~5.1
in~\cite{Thonnes}.  The second is that the poset~$\Xc$ is allowed to be ``upwardly
unbounded'' and so need not have a~$\oh$:\ 
consult, again, \cite{MS} and~\cite{Thonnes}.

(b)~Dealing with monotone rules on partially ordered state spaces without~$\oh$ is
problematic and requires the use of ``dominating processes.''
We comment that a dominating process provides a sort of {\em random\/} bounding
process and is useful when the state space is noncompact, but we shall not pursue
these ideas any further here. 
See
Kendall~\cite{Kendall} and Kendall and M{\o}ller~\cite{KM} in the context of CFTP;
we hope to discuss the use of dominating processes for our algorithm in future
work.
}
\end{remark}

\subsection{Rigorous description of Algorithm~7.2}
\label{SMrigor}

\hspace\parindent
Let~$K$ be an SM kernel with stationary distribution~$\pi$ and let~$\Kt$ be its
time-reversal, exactly as in the paragraph containing~(\ref{reversal}).  
Since we
will not be using a transition rule, we otherwise forsake the development in
Sections~\ref{app1}--\ref{imp} and apply afresh the general framework of
Section~\ref{formal}.
For simplicity, we assume $\pi(\{\zh\}) > 0$;
weaker conditions are possible for sufficiently regular
chains---see Machida~\cite{Machida2000}.
\medskip

{\em Upward kernels:\/}\ According to Theorem~1 of Kamae, Krengel, and
O'Brien~\cite{KKO}, Definition~\ref{stochdef}(b) is equivalent to the existence of
an {\em upward kernel\/}~$M$ [i.e.,\ a Markov kernel on~$\Xc$ such that, for
all~$x$, $M(x, \cdot)$ is supported on $\{y \in \Xc: y \geq x\}$] satisfying
$\nu_2 = \nu_1 M$.  Thus our assumption that~$K$ is~SM implies the existence of
upward kernels $M_{xy}$, $x \leq y$, such that
\begin{equation}
\label{up}
K(y, \cdot) = \int\,K(x, dx')\,M_{xy}(x', \cdot)\mbox{\ \ \ for all $x \leq y$.}
\end{equation}
Choose and fix such kernels $M_{xy}$; we will assume further that $M_{xy}(x', B)$
is jointly measurable in~$(x, y, x') \in \Xc^3$ for each $B \in \Bc$.
\medskip

{\em The Markov chain and the first probability space:\/}\ We consider a
probability space $(\Omega, \Ac, P)$ on which are defined a Markov chain $\vecX =
(X_0, \ldots, X_t)$ satisfying
$$
P(X_0 \in dx_0, \ldots, X_t \in dx_t) = \pi(dx_0) K(x_0, dx_1)\!\cdots\!K(x_{t-1},
dx_t)
$$
and another process~$\vecY = (Y_0, \ldots, Y_t)$ such that
\begin{equation}
\label{xy}
P(\vecY \in d\vecy\,|\,\vecX = \vecx) = \delta_{\oh}(dy_0)\,M_{x_0, \oh}(x_1,
dy_1)\,\cdots\,M_{x_{t - 1}, y_{t - 1}}(x_t, dy_t).
\end{equation}
Recalling the notational note at the end of Section~\ref{intro}, let
\begin{equation}
\label{SMCdef}
C := \{ Y_t = \zh\}
\end{equation}
and observe that conditions (i)--(iv) and~(v$'$) in Section~\ref{formal} are
satisfied by taking
\begin{eqnarray}
\label{SMdefs}
(\Xc', \Bc') &=& (\Xc, \Bc),\ \ \ B' = \{\zh\} (\in \Bc'),\ \ \ \Bc'' =
\{\emptyset, B'\}, \\
\nonumber
X &=& X_0,\ \ \ X' = X_t,\ \ \ Y = Y_t,\ \ \ \mbox{$C$ as at~(\ref{SMCdef}).}
\end{eqnarray}
Condition~(v$'$) follows from the independence of~$X_0$ and~$\vecY$, which in turn
can be verified by a simple calculation using~(\ref{xy}) and~(\ref{up}).
\medskip

{\em The second probability space and the algorithm:\/}\ To set up the second
probability space, we compute conditional probability distributions in stages, as
in Section~\ref{rigor}.  As there,
$$
P(X_0 \in dx_0, \ldots, X_{t - 1} \in dx_{t - 1}\,|\,X_t = x_t) = \Kt(x_t,
dx_{t - 1}) \cdots \Kt(x_1, dx_0)
$$
serves as a conditional distribution for $(X_0, \ldots, X_{t - 1})$ given $X_t =
x_t$.  Furthermore, (\ref{xy}) gives a conditional distribution for~$\vecY$
given~$\vecX$.

We now see how to set up and simulate our second space.  The user sets
\linebreak $\XX_t :=
\zh$, then draws an observation $\XX_{t - 1} = x_{t - 1}$ from $\Kt(x_t, \cdot)$,
then an observation $\XX_{t - 2} = x_{t - 2}$ from $\Kt(x_{t - 1}, \cdot)$, etc.\
\ Next, having chosen $\vecXX = \vecx$, the user draws an observation $\vecYY =
\vecy$ from $\Lc(\vecY\,|\,\vecX = \vecx)$.  [In detail, this is done by setting
$\YY_0 := \oh$, then drawing $\YY_1 = y_1$ from $M_{x_0, \oh}(x_1, \cdot)$, then
$\YY_2 = y_2$ from $M_{x_1, y_1}(x_2, \cdot)$, etc.]\ \ Finally, the user sets
$\CC := \{\YY_t = \zh\}$.  With the definitions~(\ref{SMdefs}), this definition
of~$\CC$, and
$$
\XX = \XX_0,\ \ \ \XXp = \XX_t,\ \ \ D(\zh) = 1 / \pi(\{\zh\}),
$$
assumptions (i)--(x) of Section~\ref{formal} are verified in a straightforward
manner.  According to Corollary~\ref{ggeneralcor}, if $\PP(\CC) > 0$, then
$\Lc(\XX_0|\CC) = \pi$.  Thus Algorithm~\ref{SMnontech} works as claimed.

\begin{remark}
\label{less}
{\em 
If the stronger assumption of realizable
monotonicity holds, then
Algorithm~\ref{SMnontech} reduces to the specialization of
Algorithm~\ref{nontech} discussed above.
This follows from the fact that one can then take
\begin{equation}
\label{RMup}
M_{xy}(x', \cdot) := P(\phi(y, U_1) \in \cdot\,|\,\phi(x, U_1) = x'),
\end{equation}
for the upward kernels in~(\ref{up}), provided that~RHS(\ref{RMup}) is jointly
measurable in $(x, y, x')$.
}
\end{remark}

\subsection{Performance of Algorithm~7.2}
\label{performance}

\hspace\parindent
Using condition~(x), we see that the routine in Algorithm~\ref{SMnontech} has
probability
\begin{equation}
\label{perform}
\PP(\CC) = \frac{P(C)}{\pi(\{\zh\})} = \frac{K^t(\oh, \{\zh\})}{\pi(\{\zh\})}
\end{equation}
of accepting the generated value~$\XX_0$.  To understand this in another way,
note that our assumption $\pi(\{\zh\}) > 0$ implies that $\Kt^t(\zh, \cdot) \ll
\pi(\cdot)$, and that the decreasing function $x \mapsto K^t(x, \{\zh\}) /
\pi(\{\zh\})$ on~$\Xc$ serves as the Radon--Nikodym derivative (RND) $x \mapsto
\Kt^t(\zh, dx) / \pi(dx)$.  With this choice of RND, we have
$$
\PP(\CC) = \inf_x \frac{\Kt^t(\zh, dx)}{\pi(dx)};
$$
this last expression is a natural (though stringent) measure of agreement between
the distribution $\Kt^t(\zh, \cdot)$ and the stationary distribution.  In the
discrete case, our performance results reduce to results found in
Sections~7--8 of Fill~\cite{Fill}; consult that paper for further discussion.

\subsection{An extension: stochastic cross-monotonicity}
\label{crosssub}

\hspace\parindent
There is no reason that the chains~$\vecXX$ and~$\vecYY$ in
Algorithm~\ref{SMnontech} need have the same kernel.  Thus, consider two (possibly
different) kernels~$K$ and~$L$ satisfying the {\em stochastic
cross-monotonicity\/} (or {\em cross-SM}) property
$$
K(x, \cdot) \preceq L(y, \cdot)\mbox{\ \ \ whenever $x \leq y$;}
$$
if $K = L$, this reduces to the Definition~\ref{stochdef}(c) of~SM.  We can then run
Algorithm~\ref{SMnontech}, replacing ``with $\YY_0 = \oh$ and kernel~$K$'' by
``with $\YY_0 = \oh$ and kernel~$L$.''  The rigorous description of this algorithm
is left to the reader.  The analogue of~(\ref{up}) is of course
\begin{equation}
\label{KLup}
L(y, \cdot) = \int\,K(x, dx')\,M_{xy}(x', \cdot)\mbox{\ \ \ for all $x \leq y$;}
\end{equation}
and if we have {\em cross-monotone transition rules\/} $\phi_K$ and~$\phi_L$
[according to the appropriate generalization of (\ref{mondef})], then one
can take
\begin{equation}
\label{cRMup}
M_{xy}(x', \cdot) := P(\phi_L(y, U_1) \in \cdot\,|\,\phi_K(x, U_1) = x').
\end{equation}
In the cross-SM case, we have the extension
$$
\PP(\CC) = \frac{P(C)}{\pi(\{\zh\})} = \frac{L^t(\oh, \{\zh\})}{\pi(\{\zh\})}
$$
of~(\ref{perform}).

One application of the cross-SM case which may arise in practice (Machida plans to
employ the following idea in a future application of our algorithm to the mixture
problems considered by Hobert et al.~\cite{HRT}) is that of {\em stochastic
dominance\/} [$K(x, \cdot) \preceq L(x, \cdot)$ for all $x \in \Xc$] by an~SM
kernel~$L$; indeed, then $K(x, \cdot) \preceq L(x, \cdot) \preceq L(y, \cdot)$
whenever $x \leq y$.  Suppose further that~$L$ has stationary
distribution~$\sigma$ with $\sigma(\{\zh\}) > 0$; then, with $\rho :=
\sigma(\{\zh\}) / \pi(\{\zh\})$,
\begin{equation}
\label{performc}
\PP(\CC) = \inf_y \frac{L^t(y, \{\zh\})}{\pi(\{\zh\})} = \rho \inf_y
\frac{\Lt^t(\{\zh\}, dy)}{\sigma(dy)},
\end{equation}
using the decreasing function $y \mapsto L^t(y, \{\zh\}) / \sigma(\{\zh\})$ as the
choice of~RND $y \mapsto \Lt^t(\{\zh\}, dy) / \sigma(dy)$.

The practical implication (say, for finite-state problems, as we shall assume
for ease of discussion for the remainder of this section) is that if the user is
unable to find a rapidly mixing~SM kernel~$K$ with stationary distribution~$\pi$
(the distribution of interest) but can find a rapidly mixing stochastically
dominant~SM kernel~$L$ with stationary distribution~$\sigma$ ``not too much
larger than~$\pi$'' (at $\{\zh\}$), then our cross-SM algorithm can be applied
efficiently, provided the imputation of~$U_1$ inherent in~(\ref{cRMup}) can be
done efficiently.

To temper enthusiasm, however, we note that the limit (namely, $\rho$) as
$t \to \infty$ \linebreak of the acceptance probability~(\ref{performc}) can often
be achieved in simpler fashion.  Suppose, for example, that simulation
from~$\sigma$ is easy, that~$x \mapsto \pi(\{x\})$ and \linebreak $x \mapsto
\sigma(\{x\})$ can be computed exactly except for normalizing constants, and that
$x = \zh$ minimizes the ratio $\sigma(\{x\}) / \pi(\{x\})$.  Then one can employ
elementary rejection sampling to simulate from~$\pi$, with acceptance
probability~$\rho$.

\sect{Relation to CFTP}
\label{taleof2}

\subsection{Comparison}
\label{comparison}
\hspace\parindent
How does our extension
of Fill's algorithm, as given by Algorithm~\ref{nontech},
compare to CFTP?  As we see it, our
algorithm has two main advantages and one main disadvantage.
\medskip

{\em Advantages:\/}\ As discussed in Section~\ref{intro} and
Remark~\ref{genremark}(b) and in~\cite{Fill},
a primary
advantage of our
Algorithms~\ref{nontech} and~\ref{altalg}
is interruptibility.  A related second advantage of Algorithm~\ref{nontech} concerns
memory allocation.  Suppose, for example, that our state space~$\Xc$ is finite and
that each time-step of Algorithm~\ref{nontech}, including the necessary imputation
(recall Section~\ref{imp}), can be carried out using a bounded amount of memory. 
Then, for fixed~$t$, our algorithm can be carried out using a fixed finite amount
of memory.  Unfortunately, it is rare in practice that the kernel~$K$ employed is
sufficiently well analyzed that one knows in advance a value of~$t$ (and a value
of the seed~$\XX_t$) giving a reasonably large probability~$\PP(\CC)$ of
acceptance.  Furthermore, the fixed amount of memory needed is in practice larger
than the typical amount of memory allocated dynamically in a run of CFTP.
See also the
discussion of read-once CFTP in Section~\ref{conn}.
\medskip

{\em Disadvantage:\/}\ A major disadvantage of our algorithms
concerns computational complexity.  We refer the reader to~\cite{Fill}
and~\cite{MTFalgo} for a more detailed discussion in the setting of realizable
monotonicity (and, more generally, of stochastic monotonicity).
Briefly, if
no attention is paid to memory usage, our algorithms have running time competitive
with CFTP: cf.~Remark~\ref{genremark}(c), and also the discussion in Remark~9.3(e)
of~\cite{Fill} that the running time of our Algorithm~\ref{nontech} is, in a certain sense,
best possible in the stochastically monotone setting.  However, this analysis
assumes that running time is measured in Markov chain steps;  unfortunately,
time-reversed steps can sometimes take longer than do forward steps to execute
(e.g.,~\cite{MTFalgo}), and the imputation described in Section~\ref{imp} is
sometimes difficult to carry out.  Moreover, the memory usage for naive
implementation of our algorithm can be exorbitant; how to trade off speed for
reduction in storage needs is described in~\cite{Fill}.

\subsection{Connection with CFTP}
\label{conn}

\hspace\parindent
There is a simple connection between CFTP and our
Algorithm~\ref{altalg}.  Indeed, suppose we carry out the usual CFTP algorithm to
sample from~$\pi$, using kernel~$K$, transition rule~$\phi$, and driving
variables~$\vecUU = (\UU_0, \UU_{-1}, \ldots)$.  Let~$\TT$ denote the backwards
coalescence time and let $\XX_0 \sim \pi$ denote the terminal state output by CFTP. 
Let $\WW \sim \pi$ independent of~$\vecUU$, and follow the trajectory from $\XX_{-\TT} :=
\WW$ to $\XX_0$; call this trajectory $\vecXX = (\XX_{-\TT}, \ldots, \XX_0)$.  Since~$\XX_0$ is
determined solely by~$\vecUU$, the random variables~$\WW$ and~$\XX_0$ are independent.

When $\pih = \pi$ in Algorithm~\ref{altalg}, the algorithm simply constructs the
same probability space as for CFTP, but with the ingredients generated in a
different chronological order:\ first $\XX_0, \XX_{-1}, \ldots$; then~$\vecUU$ (which
determines~$\TT$); then $\WW := \XX_{-\TT}$.  Again $\XX_0 \sim \pi$ and $\WW \sim \pi$ are
independent.

Note that a fundamental difference between Algorithm~\ref{altalg} and CFTP
is in what values they report.  CFTP reports $\XX_0$ as its
observation from~$\pi$,
while Algorithm~\ref{altalg} reports the value $\WW = \XX_{-\TT}$.

\begin{remark}
\label{connremark}
{\em (a)~Because of this statistical independence, it does not matter in
Algorithm~\ref{altalg} that we actually use $\XX_0 \sim \pih \neq \pi$.

(b)~The fact (1)~that $\WW$, unlike~$\XX_0$, is independent of~$\vecUU$, together with
(2)~that $\TT$ depends solely on~$\vecUU$, explains why our algorithm is interruptible
and CFTP is not.

(c)~In a single run of CFTP, the user would of course be unable to choose $\WW
\sim \pi$ as above, just as in a single run of Algorithm~\ref{altalg} we do not
actually choose $\XX_0 \sim \pi$.  So one might regard our described connection
between the two algorithms as a bit metaphorical.  But see Section~\ref{efficient}.
}
\end{remark}

{\em Read-once CFTP:\/}\
We note also that our Algorithm~\ref{nontech} bears close resemblance
to the read-once CFTP algorithm of Wilson~\cite{WilsonReadOnce}.  
Indeed, that algorithm shares
with ours the property of requiring only a bounded amount of memory.
Furthermore, in the execution of Algorithm~\ref{nontech},
consider relabeling time as follows.  At the $k^{\rm th}$
attempt at coalescence ($k = 1, 2, \ldots$),
consider $\XX_{- (k - 1) t}, \ldots, \XX_{- k t}$
in place of $\XX_t, \ldots, \XX_0$, respectively,
where now
$\XX_{-kt}$ plays two roles: the candidate for output in attempt $k$ and
the starting state in attempt $k+1$.  Then  Algorithm~\ref{nontech}
will output $\XX_{-kt}$ if and only if coalescence occurs
in attempt $k$.  Read-once CFTP, on the other hand, runs chains forwards in time
from time~$0$,
outputting $\XX_{k t}$ if and only if coalescence occurs
beginning at time~$k t$ and ending at time~$(k + 1) t$
and has also occurred on some earlier such interval.
Thus, in some sense, the two algorithms are mirrors of one another.  However, 
the read-once algorithm is {\em not} interruptible, because its output (which 
necessarily follows a previous coalescence), conditional
on rapid termination, is biased towards values resulting from fast first coalescence.
By contrast, the output of Algorithm~\ref{nontech} does not require computation of a 
previous coalescence.

\sect{Perfect samples of arbitrary size}
\label{n}

\hspace{\parindent}
Thus far in this paper we have considered only the problem of obtaining a single
observation from the distribution~$\pi$ of interest.  What can be done to obtain a
perfect sample of size~$n$?  

\subsection{Elementary options}
\label{elementary}

\hspace{\parindent}
As discussed in Section~3 of Fill~\cite{Fill} and in
Murdoch and Rosenthal~\cite{MR}, the simplest option is to apply, repeatedly and
independently, an algorithm producing a perfect sample of size~$1$.  When the
size-$1$ algorithm employed is interruptible (e.g.,\ Algorithm~\ref{nontech} or
Algorithm~\ref{altalg}), the resulting size-$n$ algorithm is both observationwise
interruptible [in the sense that, for $k = 1, \ldots, n$, the $k$th observation
output, say~$\WW_k$, and the number of Markov chain steps required,
say~$\TT_k$, are conditionally independent given all randomness used in the
generation of $\WW_1, \ldots, \WW_{k - 1}$] and totally interruptible [in the
sense that $\TT_{+} := \sum_{k = 1}^n \TT_k$ and $\vecWW := (\WW_1, \ldots,
\WW_n)$ are independent].  A related observation is that the conditional
distribution of a fixed-duration sample given its size is that of an i.i.d.\
sample, again provided that time is measured in Markov chain steps; contrast
comment~3 concerning CFTP in Remark~5.3 of~\cite{Fill}.

Another simple option is to generate a single observation, say~$\WW_1$, from a
size-$1$ perfect sampler and then run the chain~$K$ (or~$\Kt$, or any other
chain with stationary distribution~$\pi$) for $n - 1$ steps from~$\WW_1$, obtaining
$\vecWW = (\WW_1, \ldots, \WW_n)$.  Note that each observation~$\WW_k$ is
marginally distributed as~$\pi$.  This size-$n$ algorithm is also observationwise
interruptible and totally interruptible.  Of course, statistical use of~$\vecWW$
is complicated by the serial dependence of its entries.

A third option is to compromise between the first two ideas and generate $\nu$
independent vectors $\vecWW^{(1)}, \ldots, \vecWW^{(\nu)}$, where $\vecWW^{(i)}$ is
obtained using the size-$t_i$ sampler described in the preceding paragraph and
$\nu$ and $(t_1, \ldots, t_{\nu})$ are chosen (in advance, deterministically) so
that $t_1 + \cdots + t_{\nu} = n$.  Again we have both forms of interruptibility. 
This third option is an interruptible analogue of the ``RCFTP (Repeated CFTP)
tours'' of Murdoch and Rosenthal~\cite{MR}.

\subsection{Efficient use of size-$1$ perfect samplers}
\label{efficient}

\hspace\parindent
All three of the options in Section~\ref{elementary} seem wasteful.  After all, a
great deal of randomness and computational effort goes into the use of
Algorithm~\ref{nontech} or~\ref{altalg}, yet only a single observation from~$\pi$
comes out.  Is there a way to be more efficient?  Our short answer is this: Yes,
but only for Algorithm~\ref{altalg}, and then the advantage of interruptibility is
forfeited.  The remainder of this subsection provides an explanation.
\medskip

{\em Using Algorithm~{\rm \ref{nontech}}:\/}\ Suppose that we have in hand an
observation
$\WW_1 \sim \pi$, for example from a first run of Algorithm~\ref{nontech}.  If we
feed~$\WW_1$ into the routine of Algorithm~\ref{nontech} as~$\XX_t$, then the
resulting~$\vecXX$ is unconditionally distributed as a stationary trajectory
from~$K$; in particular, $\XX_s \sim \pi$ for $0 \leq s \leq t$.  Unfortunately,
this is not generally true conditionally given success (i.e.,\ given coalescence):

\begin{example}
\label{toybad}
{\em Consider again the toy random-walk example of Section~\ref{toy}, again
with~$t = 2$ but now with the modifications
\begin{eqnarray*}
&& k(0, 0) = k(2, 2) = 3/4,\ \ k(0, 1) = k(2, 1) = 1/4,\ \ k(1, 0) = k(1, 2) = 1/2,
\\
&& k(0, 2) = k(1, 1) = k(2, 0) = 0
\end{eqnarray*}
and $\pi = (2/5, 1/5, 2/5)$ to~$K$ and~$\pi$.  Suppose that the seed
value~$\XX_2$ is chosen according to~$\pi$, and that the independent-transitions
rule is used.  Then of course $\Lc(\XX_0|\CC) = \pi$, but one can check
that $\Lc(\XX_1|\CC) = (7/22, 8/22, 7/22)$ and that
$\Lc(\XX_2|\CC) = (4/11, 3/11, 4/11)$.
}
\end{example}

This example illustrates a rather catastrophic fact about Algorithm~\ref{nontech}:
If we use the output $\WW_2 := \XX_0$ as an observation from~$\pi$, we may no
longer use $\WW_1 = \XX_t$ as such.  We know of no systematic way to make use of
the auxiliary randomness generated in a run of Algorithm~\ref{nontech}.
\medskip

{\em Using Algorithm~{\rm \ref{altalg}}:\/}\ Let Algorithm~\ref{altalg}$'$
denote Algorithm~\ref{altalg} modified as follows.  For fixed~$t_0$,
use~$\TT' := \max(\TT,t_0-1)$ in place of $\TT$.  Algorithm~\ref{altalg}$'$
uses a conservative detection rule (see section~\ref{detection}), and
so is still valid.

Suppose again that we have in hand an observation distributed as~$\pi$, say from
a first run of Algorithm~\ref{altalg}$'$ or~\ref{altalg} or~\ref{nontech}.  If we
call this observation~$\WW^{(1)}_{t_0}$ and feed it into Algorithm~\ref{altalg}$'$
as~$\XX_0$, then the resulting trajectory $(\ldots, \XX_{-2}, \XX_{-1}, \XX_0)$ is
distributed as a stationary trajectory from~$K$; in particular,
$$
\vecWW^{(1)} = (\WW^{(1)}_1, \ldots, \WW^{(1)}_{t_0}) := (\XX_{-(t_0 - 1)},
\ldots, \XX_0)
$$
is a stationary trajectory of length~$t_0$ from~$K$.  Defining $\WW^{(2)}_{t_0}
:= \XX_{-\TT}$, we now feed~$\WW^{(2)}_{t_0}$ into Algorithm~\ref{altalg}$'$ as
the new seed~$\XX_0$; using randomness otherwise independent of the first run, we
obtain another stationary trajectory, call it~$\vecWW^{(2)}$, of length~$t_0$
from~$K$.  Taking~$\WW^{(3)}_{t_0}$ to be the value of~$\XX_{-\TT}$ from this
second run, we then provide~$\WW^{(3)}_{t_0}$ as a seed producing~$\vecWW^{(3)}$,
and so on.

It is not hard to see that, for
any $\nu \geq 1$, the joint distribution of the ``tours'' (as Murdoch and
Rosenthal~\cite{MR} call them) $\vecWW^{(\nu)}, \ldots, \vecWW^{(2)}, \vecWW^{(1)}$
is the same as the joint distribution of~$\nu$ serially generated ``Guarantee Time
CFTP'' (GTCFTP) tours as in Section~5 of~\cite{MR}.  (Their ``guarantee
time''~$T_g$ is our~$t_0 - 1$.)  Since each tour~$\vecWW^{(i)}$ has the marginal
distribution
\begin{equation}
\label{tourstat}
\PP(\vecWW^{(i)} \in d\vecw) = \pi(dw_1)\,K(w_1, dw_2)\,\cdots\,K(w_{t_0 - 1},
dw_{t_0})
\end{equation}
and since both $(\vecWW^{(1)}, \vecWW^{(2)}, \ldots, \vecWW^{(\nu)})$ and the
sequence of GTCFTP tours are clearly tour-valued Markov chains, we see that our
tour-chain is simply the stationary time-reversal of the (stationary) GTCFTP
tour-chain, where the common stationary tour-distribution is~(\ref{tourstat}).

Since GTCFTP tours are analyzed by Murdoch and Rosenthal~\cite{MR}, we refer the
interested reader to~\cite{MR} for further discussion and analysis.  One highlight
is that while the tours $\vecWW^{(1)}, \vecWW^{(2)}, \ldots$ are not independent,
they are $1$-dependent:\ that is, $\vecWW^{(i)}$ and $(\vecWW^{(1)},
\vecWW^{(2)}, \ldots, \vecWW^{(i - 2)}, \vecWW^{(i + 2)}, \vecWW^{(i + 3)},
\ldots)$ are independent for every~$i$.

\begin{remark}
\label{notinterruptible}
{\em For $t_0 = 1$, one can check that the tour-algorithm we have described is an
observationwise interruptible and totally interruptible algorithm for producing
{\em independent\/} observations from~$\pi$.  Unfortunately, for $t_0 \geq 2$, one
can check that all interruptibility is lost, due to the dependence of $(\XX_{-(t_0
- 1)}, \ldots, \XX_{-1})$ and~$\XX_{-T}$ (unlike the independence of~$\XX_0$
and~$\XX_{-T}$) in Algorithm~\ref{altalg}$'$.
}
\end{remark}

In light of the above remark, and the fact that CFTP tours are generated
considerably more easily than are our Algorithm~\ref{altalg}$'$ tours, we see
nothing to recommend the use of Algorithm~\ref{altalg}$'$ in practice.  
(We remark, nonetheless, that it was
in investigating our tours that we discovered the connection between
Algorithm~\ref{altalg} and CFTP described in Section~\ref{conn}.)

\bigskip\bigskip\noindent
\bf Acknowledgments.
\rm We thank the anonymous referees for many helpful suggestions about the
content and exposition of this paper, including the 
connections between Algorithm~\ref{nontech} and read-once CFTP.
We also thank David Wilson for helpful comments.

\reversemarginpar

\end{document}